# Evacuation Route Planning for Alternative Fuel Vehicles

Denissa Sari Darmawi Purba[1], Eleftheria Kontou[1*], Chrysafis Vogiatzis[2]

[1] Civil and Environmental Engineering, University of Illinois at Urbana-Champaign
[2] Industrial and Enterprise Systems Engineering, University of Illinois at Urbana-Champaign
[*] Corresponding author, email: kontou@illinois.edu, ORCiD: 0000-0003-1367-4226

**ABSTRACT**

As the number of adopted alternative fuel vehicles increases, it is crucial for communities (especially those that are susceptible to hazards) to make evacuation plans that account for such vehicles refueling needs. During emergencies that require preemptive evacuation planning, travelers using alternative fuel vehicles are vulnerable when evacuation routes do not provide access to refueling stations on their way to shelters. In this paper, we formulate and solve a novel seamless evacuation route plan problem, by designing $|K|$ spanning trees with side constraints that capture the refueling needs of each $k \in K$ vehicle fuel type on their way to reach a shelter. We develop a branch-and-price matheuristic algorithm based on column generation to solve the evacuation route planning problem. We apply the proposed framework to the Sioux Falls transportation network with existing infrastructure deployment and present numerical experiments. Specifically, we discuss the optimal system evacuation travel and refueling times under scenarios of various alternative fuel vehicle driving ranges. Our findings show that the characteristics of each vehicle fuel type, like driving range and the refueling infrastructure topology, play a pivotal role in determining evacuation route plans. This means that an evacuation route could prove unique to a single vehicle fuel type, while being infeasible to others. Finally, we observe that the driving range constraints could force evacuee vehicles to detour to meet their refueling needs before reaching safety and increase the total evacuation time by 7.32% in one of evaluated scenarios.

**INTRODUCTION**

Evacuation planning is crucial to achieving resilient and safe communities in the face of natural and anthropogenic hazards (Federal Highway Administration, 2006). Since 2018, the United States Federal Emergency Management Agency (FEMA) has announced 282 major disaster declarations, due to floods, hurricanes, and wildfires, resulting in 182 evacuation orders nationwide (FEMA, 2022). In recent years, we observed more frequent storm and hurricane emergency declarations, affecting more than 1,000 residents who need to evacuate during the Atlantic hurricane season in the southern US (National Oceanic and Atmospheric Administration, 2020). One lesson learned from Hurricane Katrina is that evacuation route planning is streamlining the evacuation process and can result in avoiding delays in emergency resources distribution (Apte, 2009; Townsend, 2006). Emergency authorities at the federal, state, and local levels are required to provide a designated evacuation route plan, leading evacuees to shelters during preemptive evacuations (Federal Highway Administration, 2017). Responding to the recent ambitious federal and state programs for accelerating alternative fuel vehicle adoption (State of Illinois, 2021; The White House, 2021), it is essential to consider these vehicles' travelers when designing evacuation routes to shelters.

Alternative fuel vehicles adoption enables the transition to decarbonized transportation systems (Li et al., 2015; Miotti et al., 2016; US Energy Information Administration, 2020). Battery electric and fuel cell electric vehicles provide competitive economic, energy security, and environmental benefits (He et al., 2019; Jacobson, 2009; Tessum et al., 2014), and serve as solutions to curb the light-duty transportation sector's dependency on fossil-fuels. Over the last decade (2011-2021), the total US sales of battery electric and fuel cell electric vehicles have reached 1.92 million (California Energy Commission, 2021). Several states in the US set aspirational goals to reduce tailpipe emissions and improve the passenger fleet's fuel efficiency. In 2018, California set out to have five million electric vehicles by 2030 and all zero-emission vehicle sales by 2050 (Lutsey, 2018). New York also committed $701 million to achieve more than 1.5 million electric vehicles and clean transit under the Climate Leadership and Community Protection Act by 2050 (NYC Office of Climate and Sustainability, 2021; Roberts, 2019).

Despite sustainability benefits, the alternative fuel vehicles market growth is subject to vehicle technology and infrastructure availability barriers. Alternative fuel vehicles have a limited driving range, depending on their fuel tank or battery pack size, sparse refueling and charging infrastructure networks, and long refueling and charging times. The driving range of an electric vehicle varies from 58 to 335 miles on a single full battery charge, and that of a fuel cell electric vehicle between 256 to 366 miles. These driving ranges fall short compared to gasoline vehicles, whose median range is 418 miles (Office of Energy Efficiency & Renewable Energy, 2018). In addition, the refueling networks for these emerging vehicle technologies remain sparse, including only 50 hydrogen, 1,635 natural gas, and 49,123 electric vehicle charging stations, compared to the 142,000 registered gasoline stations in the US, as of 2021 (US Department of Energy, 2021a). Hydrogen fuel cell and natural gas vehicles can refuel in less than 4 minutes for a 20-gallon-equivalent tank size, similar to a gasoline vehicle (US Department of Energy, 2021b, 2021c). However, charging is time-consuming for battery electric vehicle drivers, since the charging rate varies from 5 miles/hour to 240 miles/hour, depending on the charging efficiency and power level (US Department of Energy, 2021d).

These vehicle constraints may lead to alternative fuel vehicle operation challenges during natural and anthropogenic hazardous events that require preemptive evacuations (Adderly et al., 2018; Feng et al., 2020). Evacuation distances to safety often exceed the driving range of alternative fuel vehicles and frequent refueling stops may be needed. For example, during Hurricanes Katrina and Rita, a large-scale and long-distance mass displacement from Louisiana to nearby states, such as Alabama and Georgia, was observed, covering distances that exceeded 388 miles (Hori et al., 2009; Lindell et al., 2019). Metaxa-Kakavouli et al. (2018) provide evidence that the displacement of Miami residents, two days prior to Hurricane Irma, exceeded 348 miles with evacuees flow heading from Miami to Jacksonville, FL. Due to the advent of diverse vehicle fuel types with refueling accessibility limitations and diverse range constraints, such vehicle evacuees are susceptible to being stranded during emergencies, without sufficient fuel or energy to reach designated shelters (FEMA, 2019).

Emergency management authorities need to provide designated paths to alternative fuel vehicle evacuees to enable safe routing alongside refueling. The Federal Emergency Management Agency (2019) and Federal Highway Administration (2020) suggest planning evacuation routes through emergency respite sites, that may have installed gasoline stations, without specific provisions for alternative fuel vehicles.

Access to refueling stations is often neglected in evacuation route modeling literature (Gao et al., 2010). Research in this field models preemptive evacuations with no specific vehicle range and refueling considerations, where every evacuee can always reach their destination (Achrekar and Vogiatzis, 2018; Campos et al., 2012; Hasan and Van Hentenryck, 2020). However, these assumptions are not applicable to alternative fuel vehicle evacuees. Evacuation planning and coordination officials must incorporate the vehicle range and refueling dependency constraints in their contingency plans and guarantee refueling station access in their designated evacuation paths.

Since the introduction of alternative fuel vehicles in the US market, there have been no adjustments on evacuation policy due to the emerging alternative fuel vehicle technologies. Densifying the sparse alternative refueling station network would address evacuation refueling needs (Adderly et al., 2018), alleviate short driving range constraints, and facilitate more energy-efficient travel (Kontou et al., 2017, 2015); but it requires significant investments (Gnann et al., 2018). In addition, the refueling infrastructure deployment is often based on habitual operations with objectives like the maximization of profits for the network provider, not emergency operations (Ghamami et al., 2020; He et al., 2015). Even though, future investments in alternative refueling infrastructure networks are expected to increase (e.g., under the Electric Vehicle Charging Action Plan (2021)), this paper addresses an immediate need for evacuation route planning for battery electric vehicles and other competing and heterogeneous alternative fuel vehicle types.

## Our Contributions

This paper contributes to and expands the evacuation network modeling literature by developing a novel mathematical model for the evacuation route plans of a set of alternative fuel vehicles $K$, each with its own unique refueling or recharging infrastructure topology on the transportation network. The evacuation routes are designed to meet evacuation attributes, as well as the provision of reliable access to refueling infrastructure while considering that multiple different evacuation routes of alternative fuel vehicle types need to be followed simultaneously. Our research aims at answering the following critical questions for emergency agencies' needs:

1) What are the optimal evacuation route plans followed simultaneously for heterogeneous alternative fuel vehicles?
2) How does the evacuation route designation for alternative fuel vehicles differ from the one for gasoline vehicles?
3) What is the impact of the vehicle characteristics (i.e., driving range and refueling station density) on the evacuation routes planning?
4) How could emergency planners efficiently address the needs of heterogeneous vehicle fuel types during evacuation route planning?

We initiate the study by reviewing pertinent literature on evacuation planning and evacuees routing, as well as alternative vehicles' refueling modeling. We propose a new formulation of an evacuation route plan problem for alternative fuel vehicles. It determines $|K|$-minimum spanning trees (every tree rooted at the shelter) with hop constraints that model the driving range and refueling needs of each vehicle fuel type $k$ on their way to shelter. We also develop a path-based reformulation to solve the problem using a matheuristic branch-and-price algorithm with a column generation approach. We apply the proposed framework to the Sioux Falls transportation network to (i) plan the optimal evacuation routes under various scenarios of alternative fuel vehicles' driving range and refueling station deployment, and (ii) evaluate the role of alternative fuel vehicles and refueling network parameters on the evacuation performance. Finally, we use the insight gained from numerical experiments to discuss potential policy recommendations regarding the evacuation route design for alternative fuel vehicles. To our knowledge, this is the first paper to meet this timely research objective by proposing both a new mathematical formulation that describes this problem, its algorithmic solution, and its policy implications.

## LITERATURE REVIEW
### Alternative Fuel Vehicles Route Planning
Alternative fuel vehicle operations are unique due to (a) frequent refueling needs, (b) a sparse refueling and recharging infrastructure network, and (c) long refueling times. Range anxiety influences the behavior of alternative fuel vehicle drivers whose comfortable driving range explains the variance in refueling and charging decisions (Franke and Krems, 2013). Electric vehicle routes can be devised under the assumption that drivers select paths based on range anxiety and generalized costs (Agrawal et al., 2016). Jiang et al. (2012) propose a distance-constrained traffic assignment problem to model the driving range limitations of electric vehicles. This work assumes that every path's length could not exceed the driving range limit and effects the feasibility of the traffic assignment. Erdogan and Miller-Hooks (2012) propose a green vehicle routing problem that introduces the driving range and refueling station dependencies in the vehicle routing problem. Their model allows planning several stops at refueling stations to eliminate the risk of running out of fuel during the vehicle routing and shows that the route feasibility depends on the refueling station configurations. He et al. (2013) develop a network equilibrium modeling framework that captures the interactions among driving range, availability of public charging, prices of electricity, and route choices of electric vehicle drivers. Wei et al. (2018) extend the traffic assignment model to couple the transportation and the power grid network. To this date, we observe that contributions on advancing the alternative fuel vehicles routing mainly focus on accommodating habitual travel demand in the literature. There are a few studies examining alternative fuel vehicles and emergency planning, especially evacuation route planning, discussed below.

As the popularity of adopted alternative fuel vehicles rises (Rezvani et al., 2015) and the frequency of hazardous events increases (Bender et al., 2010; Eshghi and Larson, 2008), evacuation route plans that support both conventional and alternative fuel vehicles are imperative. Few studies focus on coupling the topics of alternative fuel vehicles and evacuation planning. Adderly et al. (2018) provide a comprehensive review of electric vehicles and emergency policies implementation. They point out that evacuation route feasibility is one of the critical issues for electric vehicles in evacuation planning. Feng et al. (2020) conduct a feasibility assessment of electric vehicle users to evacuate during hurricanes via a case study in Florida. Both of these works (Adderly et al., 2018; Feng et al., 2020) argue that policymakers need to consider the evacuation planning of electric vehicles and suggest potential solutions of densifying charging stations, improving battery technologies, and/or adopting hybrid vehicles. These suggestions underline the vital role of supporting refueling infrastructure and increased range of vehicle technologies that would benefit the evacuation operations of alternative fuel vehicles. There is a literature gap when it comes to examining the impact of alternative fuel vehicles' range and existing refueling infrastructure on the design of preemptive evacuation routes that we aim to bridge.

**Evacuation Route Planning**

When a disaster is imminent, a plan that evacuates people from vulnerable to safe zones would empower resilient communities (Lindell et al., 2019). Evacuation models could be categorized as macroscopic and microscopic (Hamacher and Tjandra, 2002). Microscopic approaches model individual characteristics of evacuees, their interactions, and how social and environmental factors influence their movements; they are used to guide evacuation operations. In contrast, macroscopic approaches model evacuee movements as flows in a transportation network and are typically more appropriate for evacuation planning. Hasan and Van Hentenryck (2020, 2021) examine the evacuation effectiveness of microscopic and macroscopic models and show that the macroscopic approach is consistently effective and robust in solving large-scale evacuation problems. Several macroscopic evacuation route assignment models have been developed (Hamacher and Tjandra, 2001; Lindell et al., 2019), including static routing models (e.g., Bayram et al., 2015a; Chen et al., 2012; Cova and Johnson, 2003; Yamada, 1996), dynamic routing models (e.g., Achrekar and Vogiatzis, 2018a; Ogier, 1988; Opasanon and Miller-Hooks, 2010; Xie et al., 2010), and simulation-based evacuation planning (e.g., Chen et al., 2006; Chen and Zhan, 2014; Ebihara et al., 1992; Gao et al., 2010).

In the evacuation route modeling, effective evacuation plans are expected to meet several objectives and attributes, such as safety, promptness, robustness, and seamlessness (Hasan and Van Hentenryck, 2020). To achieve safe and prompt evacuations, the objective of the evacuation path planning is often the minimization of the total evacuation or network clearance time (e.g., Bayram, 2016; Ng and Waller, 2009; Sbayti and Mahmassani, 2006). The traffic assignment for evacuation can be modeled as user equilibrium or system optimum. User equilibrium models for evacuation traffic assignment are not favored, given that the assumption that evacuees perfectly know traffic conditions during rare hazard events does not hold; instead, evacuees would follow evacuation routes as directed by state officials (Lindell and Prater, 2007; Murray-Tuite and Wolshon, 2013). Ng and Waller (2010) present a bi-level model of the evacuation route and shelter location problem that assigns evacuees to the shelter with system optimum objectives in the upper level and user equilibrium objectives in the lower level. Bayram et al. (2015) propose a so-called "constrained system optimum" traffic assignment model to couple the user equilibrium and system optimum objectives in planning the fastest evacuation routes and shelter location. Their work examines the trade-off of the user equilibrium and system optimum impacts on the evacuation performance and the fairness of the route assignment.

Another practical evacuation attribute is the route seamlessness. In this work, the term seamless path denotes an evacuation plan that is intuitive and easily followed by evacuees. A lesson learned from Hurricane Katrina indicates that seamless evacuation route planning is critical for effectively reaching safety (Apte, 2009; Townsend, 2006). During Hurricane Katrina, forking (defined as diverging evacuation routes) left evacuees confused and hesitant to enter intersections, leading to major detours and delays (Hasan and Van Hentenryck, 2020; Townsend, 2006). This is the main idea behind creating an evacuation tree, where every evacuee has a single path to the shelter, eliminating forking (Andreas and Smith, 2009). Andreas and Smith (2009) define the evacuation tree as an evacuation route design with a spanning tree structure that is rooted at the terminal or safety node. This concept suggests that no flows can be split in different directions. Achrekar and Vogiatzis (2018) extend the evacuation tree concept to accommodate a budget for coordination alongside consideration for contraflow. The seamlessness evacuation attribute appears in state and local governments' evacuation route maps, with a single direction and the minimization of forking (e.g., California FPD, 2022; FDH, 2022; LCG Traffic and Transportation Department, 2022).

Additionally, congestion due to conflicts during an evacuation plan must be minimized or completely avoided. Evacuation conflicts are observed when evacuation plans allow different flow directions to use the same roads, overloading specific areas like intersections and rendering them slow and dangerous to use (Pillac et al., 2016). Cova and Johnson (2003) design an evacuation route that minimizes the number of conflicts at intersections. Campos et al. (2012) propose an evacuation route design to define independent paths from the disaster area to multiple shelters, which minimize the traffic conflicts at intersections. A commonly used practice to optimally use the available road network and reduce conflicts is contraflow, which allows for certain streets to reverse their direction, effectively increasing the capacity of the road network leading to safety (Wolshon, 2001). Contraflow has been very well-studied (e.g., Kim et al., 2008; Kim and Shekhar, 2005; Vogiatzis et al., 2013; Xie and Turnquist, 2011). Kim et al. (2008) propose a macroscopic approach for contraflow network reconfiguration, incorporating road capacity constraints and congestion factors. Xie and Turnquist (2011) conduct lane-based routing to plan a seamless evacuation that allows for lane reversals and contraflow decisions. Contraflow has also been applied in preemptive evacuation orders in past hurricane and wildfire events. The Texas Department of Transportation created contraflow plans to strategically manage the traffic conflicts and speed up the evacuation operations in interstate routes in preparation of Hurricane Ida (TxDOT, 2021). The Florida Department of Transportation applied the contraflow strategy to control the traffic flow during preemptive evacuations for the preparation of the 2021 hurricane season (FDOT, 2021).

Studies modeling evacuation routes often assume that there are no significant vehicle range and refueling dependencies (e.g., Achrekar and Vogiatzis, 2018; Campos et al., 2012; Hasan and Van Hentenryck, 2020). Access to refueling stations is often ignored in evacuation route planning, assuming that: i) every evacuee could always reach their destination without a need to refuel; or ii) the refueling infrastructure network is strategically located to allow every evacuee to access such stations. However, these assumptions are not applicable for alternative fuel vehicles since the limited and heterogeneous vehicle ranges and refueling stations would significantly undermine alternative fuel vehicles' mobility in such settings.

Evacuation studies considering refueling station dependencies have been conducted focusing on the deployment and supply planning of refueling stations for gasoline vehicles. Gao et al. (2010) develop a simulation-optimization framework to optimally deploy gasoline stations and plan supply distributions for hurricane evacuation. Sabbaghtorkan et al. (2022) developed a so-called "idealized model" to manage optimal gasoline supplies and evacuation route with fuel capacity constraints for hurricane evacuation. In this paper, we aim to develop the optimal evacuation route planning considering the vehicle range and refueling dependencies for battery electric vehicles and other competing and heterogeneous alternative fuel vehicle types.

**Bridging the Knowledge Gap**

We uncover a knowledge gap since there are no studies that couple the alternative fuel vehicles and their refueling networks specifications with evacuation route planning principles. **Table 1** highlights several studies focusing either on evacuation or alternative fuel vehicle routing. Our paper plans evacuation routes for alternative fuel vehicles and bridges this literature gap.

**Table 1 Literature Review of Evacuation and Alternative Fuel Vehicle Routing Models**

| Author | Evacuation Attributes | Alternative Fuel Vehicle Attributes |
|---|---|---|
| Jiang et al. (2012) | | Range constraints |
| Erdogan and Miller-Hooks (2012) | | Range constraints; refueling stations dependency |
| He et al. (2013) | | Range constraints; refueling stations dependency; electric charging price |
| He et al. (2013) | | Range constraints; refueling stations dependency; electric charging price; refueling time |
| Wei et al. (2018) | | Range constraints; refueling stations dependency; electric charging price; energy allocation |
| Cova and Johnson (2003) | System optimum; minimize conflict | |
| Campos et al. (2012) | System optimum; minimize conflict | |
| Bayram et al. (2015) | "Constrained System optimum" | |
| Andreas and Smith (2009) | Minimize delay risk; seamless | |
| Achrekar and Vogiatzis (2018) | Minimize delay risk; seamless; contraflow | |
| Kim et al. (2008) | System optimum; contraflow | |
| Feng et al. (2020) | User equilibrium | Range constraints; refueling station dependency |
| Sabbaghtorkan et al. (2022) | Maximize evacuee demand | Range constraints; refueling stations dependency; refueling time |
| **Our contribution** | **System optimum; seamless; minimized conflict; contraflow;** | **Range constraints; refueling stations dependency; refueling time** |

**PROBLEM FORMULATION**

We present a novel evacuation planning framework for multiple types of vehicle fuels: the $|K|$-evacuation tree routes planning problem. We aim to design evacuation paths to shelters that form evacuation trees for each $k \in K$ different vehicle fuel types. The novelty of the framework relies on the integration of the driving range constraints, vehicle refueling dependencies, and evacuation attributes when devising evacuation routes. This formulation is inspired by a well-known graph-theoretic problem that designs a centralized network with quality-of-service constraints: the minimum spanning tree problem with hop constraints (Gouveia et al., 2011, 2008).

   Our proposed evacuation model assigns an evacuation path to each evacuation origin and each vehicle fuel type in the region that carries five desirable evacuation attributes, presented below.

- **Fast evacuation**: the evacuation plan must evacuate each vehicle type's evacuees using a fast route to the shelter. System optimum network models are used for optimal evacuation paths design, where the objective is to minimize the total system's evacuation time (Bayram et al., 2015; Lindell et al., 2019).

- **Seamless evacuation route**: this guarantees that each evacuee vehicle fuel type follows one route to safety, streamlining coordination and enforcement efforts. We achieve seamlessness by making sure that the evacuation route forms a tree (Achrekar and Vogiatzis, 2018; Andreas and Smith, 2009).
- **Contraflow (lane reversal) traffic control**: the optimal use of the available road network is enforced by allowing certain streets to reverse their direction and effectively increase the network's capacity to lead to safety (Wolshon, 2001).
- **Guarantee access to a refueling station**: each evacuation path must ensure access to refueling stations to fulfill refueling needs due to alternative fuel vehicles' short driving ranges. Our study considers the driving range variation based on the vehicle fuel type (Kuby and Lim, 2005) and uses hops as proxy to represent distance measurement in the resultant evacuation network (Gouveia et al., 2011).
- **Simultaneous evacuation route plan**: we consider simultaneous evacuation routing for multiple vehicle fuel types with different driving ranges and refueling stations based on their fuel type. Our proposed formulation will determine $|K|$ evacuation trees, one for each alternative fuel vehicle type $k \in K$, followed simultaneously.

In our study, we assume that every evacuee starts evacuating with a fully fueled vehicle. This assumption is applicable to a preemptive evacuation order, given that evacuees are well-informed regarding the order. Evacuees follow the directed evacuation routes, as these are released and mandated by emergency management and coordination officials. In the remainder of this section, we define the notation, provide a formal problem statement, and present the mathematical formulation of the problem.

**Definitions and notation**

Let $G = (N, A)$ be a transportation network, where $N$ represents the set of nodes and $A$ the set of links. A safety node on the transportation network is represented by $s$. We define a set of vehicles $K$, representing each of the different vehicle fuel types. As an example, we could consider $K = \{$gasoline, hydrogen, electric$\}$ representing $|K| = 3$ different fuel vehicle types.

For every node $i \in N$ in the network, we assume that $q_i^k$ is the number of vehicles of type $k$ originating from node $i$ that need to be routed to the safety node $s$. A subset of network nodes may host refueling infrastructure. We define $ST(k) \subseteq N, \forall\, k \in K$ as the set of all nodes with installed infrastructure to refuel vehicles of type $k$. Moreover, let $z_i^k$ be a binary indicator parameter such that when node $i \in ST(k)$ (i.e., node $i$ has the necessary refueling station for vehicle type $k$), then $z_i^k = 1$; otherwise, $z_i^k = 0$. For every link $(i, j) \in A$ in the network, the travel time is denoted as $t_{ij}$ and is computed by the Bureau of Public Road (BPR) function, associating the travel time of a link with its flow to capacity ($u_{ij}$) ratio (Sheffi, 1985). Wherever necessary, we use $M$ as a big-M value in the formulation.

**Figure 1** illustrates the parameters in our problem's formulation. We enforce evacuation routing with a tree network structure; hence, all routes will form trees rooted at a single safety node. In real-world evacuation cases, we may have multiple shelters. Our model formulation can handle multiple safety nodes by creating a dummy node ("super safety" node) that serves as $s$. In that case, each existing safety node will be directly connected to $s$ with an arc of zero travel time and infinite capacity.

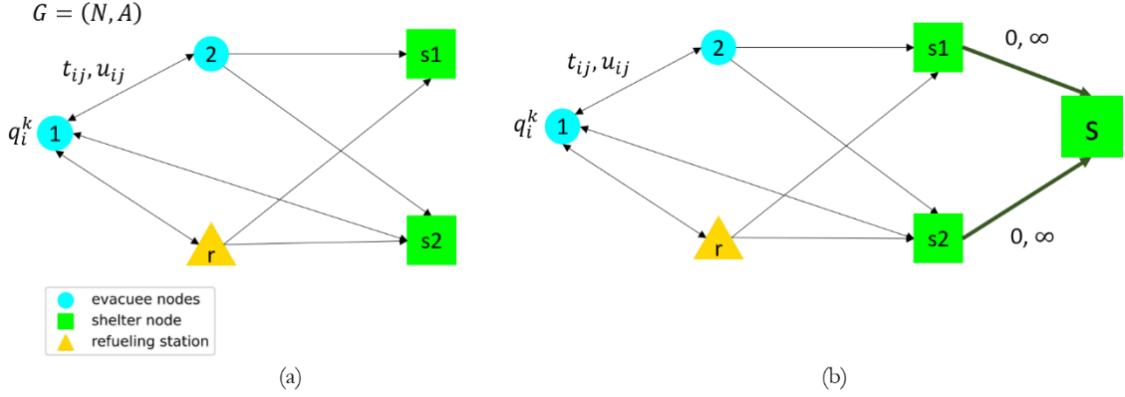

**Figure 1** Example of (a) network topology with existing refueling station and multiple shelter nodes and (b) modified network topology with a super safety node.

In the proposed evacuation problem, we model the driving range constraint for each vehicle fuel type $k$. We achieve that by imposing hop constraints that approximate the distance or length of a route. We are inspired by the concept of hop constraints that is often used to control the quality-of-service in communications network design (Gouveia et al., 2011, 2008). Hop constraints use the number of links/hops as a generic controlled parameter to manage the network topology in a macroscopic level. In the evacuation setting, the hop constraints would be beneficial in two ways: (i) they provide emergency planners a way to integrate heterogeneous fuel types in their plans; and (ii) they account for variability in driving ranges for different vehicle manufacturers.

The number of hops is denoted by $\tau_k$. It is a natural number of the transportation network links that a vehicle needs to traverse; as such it serves as a proxy to the alternative fuel vehicle's driving range. Specifically, a vehicle of fuel type $k$ will need to stop and refuel on its way to a safety node if it originates from a node that is situated at least $\tau_k$ hops away in the current evacuation plan. Additionally, evacuation planners could define the distance of one hop to correspond to the scale of distance used in the evacuation planning decisions.

Furthermore, we assume that vehicles of fuel type $k$ that need to stop at some refueling node $i \in ST(k)$ will spend refueling time proportional to the number of hops that the vehicle needs to traverse prior to arriving to safety. Refueling at a refueling station type $k$ would result in total refueling time which is a linear function of a fixed refueling efficiency rate, $r_c^k$. We assume that every station has enough fuel and parking capacity to serve all necessary vehicles during the evacuation process without congestion or queues formed on the network.

There are five sets of decision variables in this model. The first three are associated with the vehicle flows. The binary variable $x_{ij}^k$ denotes whether a link is part of the evacuation tree route for vehicle type $k$. If a link is selected to be in the evacuation tree for a certain vehicle type $k$, this implies that it would accommodate evacuation flow of evacuees with vehicle type $k$, defined as $f_{ij}^k$. The total vehicular flow traversing a link is defined as $v_{ij}$.

The remaining decision variables are associated with the refueling requirements. First, we define binary variables $w_i^k$ to denote whether vehicles of fuel type $k$ originating from node $i$ have to travel through refueling infrastructure before reaching safety. Specifically, we let $w_i^k = 1$ signal two conditions: (i) the evacuees using vehicles of type $k$ departing from node $i \in N$ need to refuel because of their starting location being more than $\tau_k$ hops away from the shelter; (ii) the evacuees using vehicles of fuel type $k$ departing from node $i \in N$ shall reroute to refueling stations due to the tree structure of the evacuation network. Condition (ii) is meant to accommodate the refueling requirement of upstream drivers that are evacuated through node $i \in N$ having originated from a node farther than $\tau_k$ hops away from the shelter. To facilitate this, we define a second binary variable $y_{i,l}^k$, which is equal to 1 if and only if node $i \in N$ is located $l$ hops away from $s$ for the evacuation tree of $k$.

**Figure 2** presents an example of the refueling and routing decisions in our proposed problem setup. Consider a network topology with one refueling station located at node $r$ and one shelter located at node $s$. In this example, we consider one vehicle fuel type with a driving range represented as $\tau_k = 3$ hops. Given the network topology, every vehicle originating from node $i$ (located more than $\tau_k$ hops away) needs to go through refueling station $r$ before reaching the shelter. Hence, we have that $w_i^k = 1$. For node $j$, there exists a path to reach $s$ whose distance is less than $\tau_k$ hops away (through node $m$). However, our model guarantees a seamless evacuation route and shall maintain an evacuation tree structure. Due to the given network topology, all vehicles from node $i$ will pass through node $j$ to access the refueling station $r$. Maintaining the seamlessness attribute, all vehicles from node $j$ must be rerouted through the refueling station to provide refueling access for its preceding evacuee nodes. Hence, we also have that $w_j^k = 1$. The same is not true for node $m$, whose distance is less than $\tau_k$ hops from the shelter. Additionally, no other vehicles that would need to access a refueling station go through $m$. Hence, we have that $w_m^k = 0$.

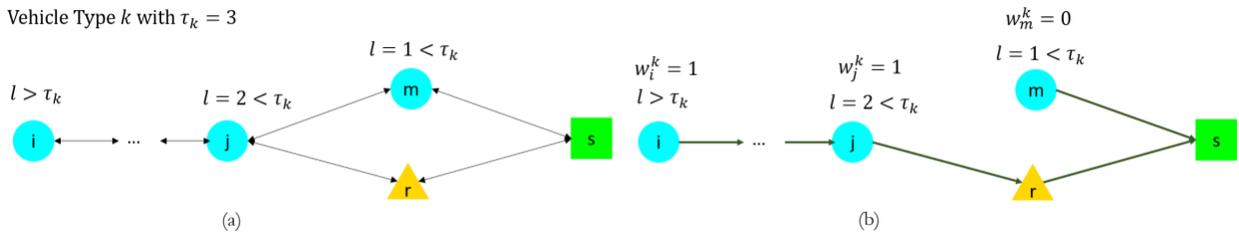

**Figure 2 The example of (a) a transportation network for vehicle type $k$ and driving range as $\tau_k = 3$ hops; (b) the refueling routing decision and the optimal evacuation routing result. Evacuee flows from node $i$ need to refuel due to their driving ranges' constraint, while evacuees from node $j$ are required to pass through a refueling station to maintain the evacuation tree structure.**

## Mathematical Formulation

We provide the formulation of the problem as in (1)—(18).

$$\min \sum_{(i,j)\in A} t_{ij}(v_{ij})v_{ij} + \sum_{k\in K}\sum_{i\in N} q_i^k r_c^k \sum_{l\geq\tau_k} y_{i,l}^k l$$

$$= \min \sum_{(i,j)\in A} t_{0,ij}\left(1+\alpha\left(\frac{v_{ij}}{u_{ij}}\right)^\beta\right)v_{ij} + \sum_{k\in K}\sum_{i\in N} q_i^k r_c^k \sum_{l\geq\tau_k} y_{i,l}^k l \tag{1}$$

$s.t.$

$$v_{ij} = \sum_{k\in K} f_{ij}^k \qquad\qquad \forall\,(i,j)\in A \tag{2}$$

$$\sum_{j:(i,j)\in A} f_{ij}^k - \sum_{j:(j,i)\in A} f_{ji}^k = q_i^k \qquad\qquad \forall\,i\in N\setminus\{s\}, \forall k\in K \tag{3}$$

$$f_{ij}^k \leq M\,x_{ij}^k \qquad\qquad \forall(i,j)\in A, \forall k\in K \tag{4}$$

$$\sum_{j:(i,j)\in A} x_{ij}^k = 1 \qquad\qquad \forall\,i\in N\setminus\{s\}, \forall k\in K \tag{5}$$

$$\sum_{(i,j)\in A:i,j\in S} x_{ij}^k \leq |S|-1 \qquad\qquad \forall S\in C \tag{6}$$

$$x_{ij}^k + x_{ji}^q \leq 1 \qquad\qquad \forall(i,j)\in A, \forall k,q\in K \tag{7}$$

$$\sum_{l=1}^{|N|-1} y_{i,l}^k = 1 \qquad\qquad \forall\,i\in N, \forall k\in K \tag{8}$$

$$y_{s,0}^k = 1 \qquad\qquad \forall k\in K \tag{9}$$

$$y_{i,l}^k \leq \sum_{j:(i,j)\in A} y_{j,l-1}^k x_{ij}^k \qquad\qquad \forall\,i\in N, \forall k\in K, \forall l=1,2,\dots,|N|-1 \tag{10}$$

$$w_s^k = 0 \qquad\qquad \forall k\in K \tag{11}$$

$$w_i^k \geq \sum_{l \geq \tau_k} y_{i,l}^k \qquad\qquad \forall i \in N \setminus \{s\}, \forall k \in K \quad (12)$$

$$w_i^k - z_i^k \leq \sum_{j:(i,j) \in A} x_{ij}^k(z_j^k + w_j^k) \qquad\qquad \forall i \in N, \forall k \in K \quad (13)$$

$$x_{ij}^k \in \{0,1\} \qquad\qquad \forall (i,j) \in A, \forall k \in K \quad (14)$$

$$y_{i,l}^k \in \{0,1\} \qquad\qquad \forall i \in N, \forall k \in K, \forall l = 1,2, \dots, |N| - 1 \quad (15)$$

$$w_i^k \in \{0,1\} \qquad\qquad \forall i \in N, \forall k \in K \quad (16)$$

$$v_{ij} \geq 0 \qquad\qquad \forall (i,j) \in A \quad (17)$$

$$f_{ij}^k \geq 0 \qquad\qquad \forall (i,j) \in A, \forall k \in K. \quad (18)$$

The objective function in (1) aims to minimize the total system evacuation time by summing the time spent traversing each link on the evacuation network and the time spent refueling. Constraints (2)—(4) are traffic assignment constraints that keep track of the total vehicle flow, enforce flow conservation on each node, and allow flow if and only if the corresponding link is part of an evacuation tree, respectively. Constraints in (5)—(6) are tree constraints. They restrict the number of outgoing links from each node to be equal to 1, as each node should have one path to safety. We have subtour elimination constraints for every cycle $S \in C$, where $C$ is the set of all cycles present in the transportation network. This set grows exponentially in cardinality. Constraints (7) are conflict and contraflow constraints to prohibit the same link to serve flow in both directions. We enforce contraflow constraints that take full advantage of the network links' capacities.

Constraints (8)—(10) define the distance, in number of hops, that vehicle flow from each node has to traverse prior to reaching safety in an evacuation tree. More specifically, constraints (8) enforce that there is exactly one distance to the shelter in each evacuation tree $k$ for each evacuee node. This is expected as from each node there is exactly one path to safety, due to the spanning tree network structure. In constraints (9), we enforce that the shelter node is the only one at zero hops from safety. Constraints (10) enforce that a node can only be located $l$ hops away from the shelter if there exists a node that is located before it in the evacuation tree that is located $l - 1$ hops away.

Constraints (11)—(13) control the detour decision through a refueling station before reaching safety. Constraints (11) enforce evacuees of any vehicle type located at the shelter node not to refuel since they have already reached safety. Constraints (12) set every vehicle located $\tau_k$ or more hops away from the shelter to pass through a refueling station and refuel. Constraints (13) control the route detours, due to refueling needs and the evacuation path's seamlessness, as discussed in **Figure 2**. The left-hand side of (13) states that if a node has a refueling station, all evacuees with such a vehicle fuel type traversing that node could meet their refueling requirement directly without detour. The right-hand side of (13) enforces that evacuation demand, which emerges from nodes that warrant refueling, must detour and go through a refueling station node. If there is no adjacent refueling station, then the subsequent node on the evacuation route must detour to refueling station to provide refueling access for the preceding evacuee nodes. Finally, the decision variable restrictions, as those are provided in the definitions section, are enforced in constraints (14)—(18).

## SOLUTION METHOD

In this section, we present a matheuristic method inspired by branch-and-price for solving the $|K|$-evacuation tree routes planning problem. Our problem is formulated as a mixed-integer nonlinear optimization problem, which is computationally intractable as the size of the network increases. Branch-and-price with column generation has been used to provide both heuristic and exact solutions to the (relevant to our problem) minimum spanning tree with side constraints problem (Fischetti et al., 2002; Tilk and Irnich, 2018).

Note that in our formulation, the nonlinearity exists due to the selection of the objective function, where the link performance is modeled after the BPR function. The travel time of each link would dynamically change during the traffic assignment process. Implementing branch-and-price with a traditional column generation scheme would be challenging as it only allows for a fixed cost parameter in each iteration. Thus, we design a matheuristic based on branch-and-price, to enable dynamic changes of travel time in the network as we improve the traffic assignment during each iteration.

In the remainder of this section, we first derive a path-based formulation for our problem. Then, we present the branch-and-price framework and propose a column generation inspired matheuristic to help solve our reformulation. For completeness, we also show the restricted main problem and the pricing subproblem we obtained for the column generation scheme.

**Path-based Reformulation**

One possible way to decompose the original $|K|$-evacuation tree routes planning problem is to reformulate it as an optimization problem with path-based decision variables. In the evacuation tree problem, for each origin, all evacuees are routed via a single path towards the shelter node. We could replace the link flow variable of each vehicle fuel type with a path flow variable representing the flow from each origin node to the shelter node for each vehicle fuel type $k$. Hence, the reformulation would focus on reconstructing traffic assignment constraints (2)—(4) into path-based ones, and couple the reformulation constraints with the remaining original constraints (including the tree constraints, contraflow constraints, and hop constraints) in the $|K|$-evacuation tree routing problem.

We introduce new additional parameters and variables for the reformulation, on top of the ones defined earlier. The set of evacuee origin nodes heading to shelter $s$ is defined as $O$. For every evacuee's origin node $o \in O$ in the network, we assume that $q_o^k$ is the number of vehicles of type $k$ originating from node $o$ that need to be routed to the safety node $s$. Each evacuee originating from node $o \in O$ would traverse a single path for each vehicle fuel type $k$ to the shelter. We define $P(k, o)$ as the set of paths between nodes $o$ and $s$ for an evacuee of vehicle type $k$. For each path $p \in P(k, o)$, $\delta_{ijp}^{ko}$ is a binary indicator parameter such that when arc $(i, j)$ belongs to path $p$ of origin node $o$ for fuel type $k$, then $\delta_{ijp}^{ko} = 1$; otherwise, $\delta_{ijp}^{ko} = 0$. Moreover, we define $\phi_p^{ko}$ as the path length of path $p$ for evacuee originating from node $o$ using vehicle fuel type $k$. Finally, we introduce new binary decision variables, $\lambda_p^{ko}$, which denote whether path $p$ is part of the evacuation tree route for evacuees originating from node $o$ using vehicles $k$. We present the path-based reformulation in (19)—(23).

$$\min \sum_{(i,j) \in A} t_{ij}(v_{ij}) v_{ij} + \sum_{k \in K} \sum_{o \in O} q_o^k r_c^k \sum_{l \geq \tau_k} y_{o,l}^k l$$
$$= \min \sum_{k \in K} \sum_{o \in O} \sum_{p \in P(k,o)} \sum_{(i,j) \in A} t_{ij}(v_{ij}) \delta_{ijp}^{ko} q_o^k \lambda_p^{ko} + \sum_{k \in K} \sum_{o \in O} \sum_{p \in P(k,o)} q_o^k r_c^k \phi_p^{ko} \lambda_p^{ko} \quad (19)$$

$s.t. \ (2), (4) — (18),$

$$f_{ij}^k = \sum_{o \in O} \sum_{p \in P(k,o)} \lambda_p^{ko} \delta_{ijp}^{ko} q_o^k \qquad \forall (i,j) \in A, \forall k \in K \quad (20)$$

$$\sum_{p \in P(k,o)} \lambda_p^{ko} = 1 \qquad \forall o \in O, \forall k \in K \quad (21)$$

$$\phi_p^{ko} = \sum_{l=\tau_k}^{|N|-1} y_{o,l}^k l \qquad \forall p \in P(k,o), \forall o \in O, \forall k \in K \quad (22)$$

$$\lambda_p^{ko} \in \{0,1\} \qquad \forall p \in P(k,o), \forall o \in O, \forall k \in K \quad (23)$$

The objective function in (19) aims to minimize the total evacuation time by summing the time spent traversing each link on the evacuation network and the time spent refueling in a path-based form. Constraints (20)—(21) are traffic assignment constraints that replace constraints (3). Constraints (20) keep track of the total flow of each vehicle fuel type $k \in K$ in arc $(i, j) \in A$. These constraints define the transformation function between the path-based flow variables, $\lambda_p^{ko}$, and the arc-based flow variables, $f_{ij}^k$. Constraints (21) guarantee that each evacuee originating from node $o$ using vehicle fuel type $k$ would follow one path to node $s$. Constraints (22) keep track of the evacuation distance of each path used by evacuees in the network. Finally, the binary decision variable restriction is enforced in constraints (23). The path-based reformulation is a mixed-integer nonlinear problem and is amenable to column generation, which we discuss next.

**Branch-and-Price with Column Generation Matheuristic**

In our proposed solution method, branch-and-bound handles the integrality of the original problem and column generation implementations involve finding the best evacuation tree to enter, where each path minimizes the evacuation time and satisfies the hop constraints. Our matheuristic approach is designed to handle the nonlinearity of the BPR function in the objective function. A matheuristic approach is the combination of mathematical programming and heuristics, allowing satisfactory solutions in mathematical programming (Archetti and Speranza, 2014). We adopt the restricted main problem matheuristic approach to incorporate a local search with fixed parameters for each branch-and-bound tree node (Archetti et al., 2011; Archetti and Speranza, 2014). Specifically, our column generation matheuristic uses a fixed parameter for the time in the restricted main problem and pricing subproblem to find a high-quality evacuation tree to enter. After identifying the entering columns, we apply a traditional branch-and-bound scheme so that we may explore better solutions.

*Restricted Main Problem*

In the Restricted Main Problem (RMP), we solve the path-based reformulation with a restricted set of columns. We only consider the feasible set of paths $p \in P(k, o)$ that complies with the driving range constraints of vehicle fuel type $k$. This set of feasible paths is generated for each pricing problem. The RMP is formulated using the path-based traffic assignment constraints, tree constraints, and contraflow constraints to seek the minimum evacuation tree route plan from the set of feasible paths.

We propose a matheuristic with a fixed time parameter in each iteration. We tackle the dynamic parameter of the original problem by updating the parameter values in each iteration. We define $\hat{v}_{ij}$ as the fixed value of the current total link flow in the ongoing traffic assignment process and $\hat{\phi}_p^{ko}$ as the fixed value of evacuation distance of path $p$ for evacuees from node $o$ using vehicles of type $k$ in the iteration. Given that we fixed the flow assignments in the restricted main problem, the RMP problem structure changes from mixed-integer nonlinear problem to be a mixed-integer linear problem. We then have:

$$\min \sum_{k \in K} \sum_{o \in O} \sum_{p \in P(k,o)} \sum_{(i,j) \in A} t_{ij}(\hat{v}_{ij}) \, \delta_{ijp}^{ko} q_o^k \lambda_p^{ko} + \sum_{k \in K} \sum_{o \in O} \sum_{p \in P(k,o)} q_o^k r_c^k \hat{\phi}_p^{ko} \lambda_p^{ko} \qquad (24)$$

$$s.t. \ (4) - (7), (14), (18), (20), (21), and \ (23)$$

Relaxing the integrality of $x_{ij}^k$ and $\lambda_p^{ko}$, we solve the relaxed RMP as linear problem. In this problem, $x_{ij}^k$ is not objective decision variable on objective function and it interlinks to $\lambda_p^{ko}$; thus, $x_{ij}^k(\lambda_p^{ko})$. We claim that if $\lambda_p^{ko^*}$ are optimal solution then the corresponding $x_{ij}^k(\lambda_p^{ko^*})$ is also optimal. Let $\pi_{ij}^k$ and $\mu_o^k$ be the dual variables associated with constraints (20) and (21), respectively. Then, $\lambda_p^{ko^*}$, $\pi_{ij}^{k^*}$, and $\mu_o^{k^*}$ are optimal if and only if it satisfies the three optimality conditions of the relaxed RMP problem:

- *Primal feasibility*

$$f_{ij}^k = \sum_{o \in O} \sum_{p \in P(k,o)} \lambda_p^{ko^*} \delta_{ijp}^{ko} q_o^k \qquad \qquad \forall \, (i,j) \in A, \ \forall k \in K \quad (25)$$

$$\sum_{p \in P(k,o)} \lambda_p^{ko^*} = 1 \qquad \qquad \forall \, o \in O, \ \forall k \in K \quad (26)$$

$$\lambda_p^{ko^*} \geq 0 \qquad \qquad \forall \, p \in P(k,o), \forall \, o \in O, \ \forall k \in K \quad (27)$$

- *Complementary slackness*

$$\pi_{ij}^{k^*} \left( \sum_{o \in O} \sum_{p \in P(k,o)} \lambda_p^{ko^*} \delta_{ijp}^{ko} q_o^k - f_{ij}^k \right) \qquad \qquad \forall \, (i,j) \in A, \ \forall k \in K \quad (28)$$

$$\mu_o^{k^*} \left( \sum_{p \in P(k,o)} \lambda_p^{ko^*} - 1 \right) \qquad \qquad \forall \, o \in O, \ \forall k \in K \quad (29)$$

- *Dual feasibility*

$$\left(\sum_{(i,j)\in A} t_{ij}(\hat{v}_{ij})\,\delta_{ijp}^{ko}q_o^k \;+\; q_o^k r_c^k \hat{\phi}_p^{ko}\right) - \sum_{(i,j)\in A} \pi_{ij}^{k^*}\delta_{ijp}^{ko}q_o^k - \mu_o^{k^*} \geq 0 \qquad \forall\, p \in P(k,o), \forall\, o \in O,\, \forall\, k \in K \quad (30)$$

Using the dual feasibility condition, we would like to find an improving path $p \in P(k,o)$ that improves the relaxed RMP. Hence, we construct the pricing subproblem that minimize the dual condition (30) for each origin $o \in O$ and vehicle type $k \in K$.

### Pricing Subproblem

This section describes the pricing subproblem (SP) of our matheuristic column generation. We only consider the feasible set of paths $p \in P(k,o)$ that complies with the range limitations for type $k$. We then formulate the pricing subproblem as a side-constrained shortest path problem, where hop constraints serve as the side constraints.

Let $\pi_{ij}^k$ and $\mu_o^k$ be the dual variables associated with constraints (20) and (21), respectively. We define $x_{ij}^{ko}$ as the binary decision variable indicating whether the link is part of the feasible evacuation path and $w_i^k$ as the binary decision variable indicating whether node $i$ need to provide access to refueling station of type $k$. As we model the evacuation distance in hop-distance unit, we could formulate $\phi_p^{ko} = \sum_{l=\tau_k}^{|N|-1} y_{o,l}^k = \sum_{(i,j)\in A} x_{ij}^{ko} w_o^k$. The pricing subproblem is formulated for each origin $o \in O$ and each fuel type $k \in K$, as follows:

$\forall\, o \in O, \forall\, k \in K, solve:$

$$\min \left(\sum_{(i,j)\in A} t_{ij}(\hat{v}_{ij}) x_{ij}^{ko}q_o^k \;+\; q_o^k r_c^k \hat{\phi}_p^{ko}\right) - \sum_{(i,j)\in A} \pi_{ij}^k x_{ij}^{ko}q_o^k - \mu_o^k$$

$$= \min \left(\sum_{(i,j)\in A} t_{ij}(\hat{v}_{ij}) x_{ij}^{ko}q_o^k \;+\; \sum_{(i,j)\in A} r_c^k x_{ij}^{ko} w_o^k q_o^k\right) - \sum_{(i,j)\in A} \pi_{ij}^k x_{ij}^{ko}q_o^k - \mu_o^k$$

$$= \min \sum_{(i,j)\in A} \left(t_{ij}(\hat{v}_{ij}) + r_c^k w_o^k - \pi_{ij}^k\right) x_{ij}^{ko}q_o^k \;-\; \mu_o^k \tag{31}$$

$s.t.$

$$\sum_{j:(i,j)\in A} x_{ij}^{ko} = 1 \qquad\qquad\qquad \forall i = o \tag{32}$$

$$\sum_{j:(i,j)\in A} x_{ij}^{ko} - \sum_{j:(j,i)\in A} x_{ji}^{ko} = 0 \qquad\qquad \forall\, i \in N\backslash\{o,s\} \tag{33}$$

$$\sum_{j:(j,i)\in A} x_{ji}^{ko} = -1 \qquad\qquad\qquad \forall i = \{s\} \tag{34}$$

$$w_s^k = 0 \qquad\qquad\qquad \forall i = s \tag{35}$$

$$\sum_{(i,j)\in A} x_{ij}^{ko} - \tau_k \;\leq\; |N| w_o^k \qquad\qquad \forall i = o \tag{36}$$

$$w_i^k - z_i^k \leq \sum_{j:(i,j)\in A} x_{ij}^k(z_j^k + w_j^k) \qquad\qquad \forall i \in N \tag{37}$$

$$x_{ij}^{ko} \in \{0,1\} \qquad\qquad\qquad \forall\,(i,j) \in A \tag{38}$$

$$w_i^k \in \{0,1\} \qquad\qquad\qquad \forall i \in N \tag{39}$$

The objective function in (31) tracks the reduced cost (pricing) for a new improving column, i.e., the feasible evacuation path for each evacuee demand node and vehicle fuel type. Constraints (32)—(34) generate paths originating from evacuees origin node to safety. Constraints (35) maintain that shelter node does not need to provide access to refueling station. Constraints (36) guarantee refueling decision; a vehicle located $\tau_k$ or more hops away from safety needs to refuel to the refueling station. Constraints (37) control the route detours due to refueling decisions. Constraints (38) and (39) reflect the binary nature of the decision variables.

### Matheuristic Column Generation Algorithm

We now present the column generation matheuristic scheme for solving the path-based reformulation problem. The matheuristic approach enables us to dynamically alter the travel time values in the network, after we conduct traffic assignment at each column iterations. The RMP is a mixed-integer program, hence, we need to first relax the integrality of the problem into a linear problem. The integrality would be solved using a branch-and-bound algorithm. With the respective linear relaxation problem, we then devise the column-generation algorithm as follows and illustrated in **Figure 3**.

**Step 0:** Start with $\hat{v}_{ij} = 0$ and $t_{ij}(\hat{v}_{ij}) = t_{ij}(0) \; \forall \, (i,j) \in A$. Construct one dummy column or feasible path $p_o^k \; \forall \, o \in O, \; k \in K$ with its corresponding evacuation path distance, $\hat{\phi}_p^{ko} = M$.

**Step 1:** Solve the RMP and calculate the dual prices $\pi_{ij}^k$ and $\mu_o^k$.

**Step 2:** Solve the SP $\forall \, o \in O, \; k \in K$

- If the reduced-price objective (25) $< 0$, add the new column to RMP, update $\hat{v}_{ij}$ with current flow assignment and $\hat{\phi}_p^{ko}$ with the current path distance in the RMP, and go back to Step 1.

- If the reduced-price objective (25) $\geq 0$, stop.

Note that as the number of feasible paths in a network is finite, the above procedure terminates in a finite number of steps.

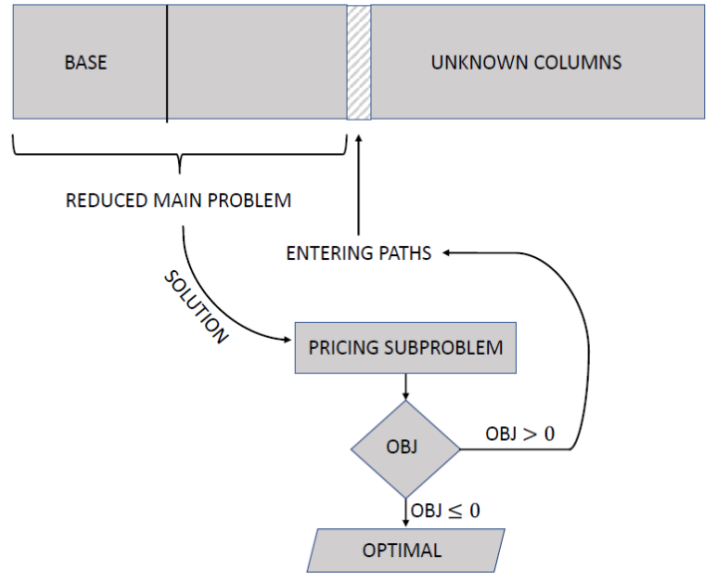

**Figure 3 Schema of the proposed column generation algorithm, adapted from** (Desrosiers and Lübbecke, 2005).

### Branch-and-Price Algorithm

In this section, we directly adopt a well-known branch-and-price generic algorithm to solve our mixed-integer problem (Desaulniers et al., 2005; Vanderbeck, 2005, 2000). As mentioned earlier, we apply branch-and-bound to handle integrality and the matheuristic column generation to solve the minimum evacuation tree iteratively. We first relax the restricted main problem to execute the matheuristic column generation. With it, we solve the relaxed RMP and evaluate the integrality of the obtained solution. If the result is not integer, we branch at the current solution and further perform the matheuristic for each branch. The branch-and-price procedure is illustrated in **Figure 4**.

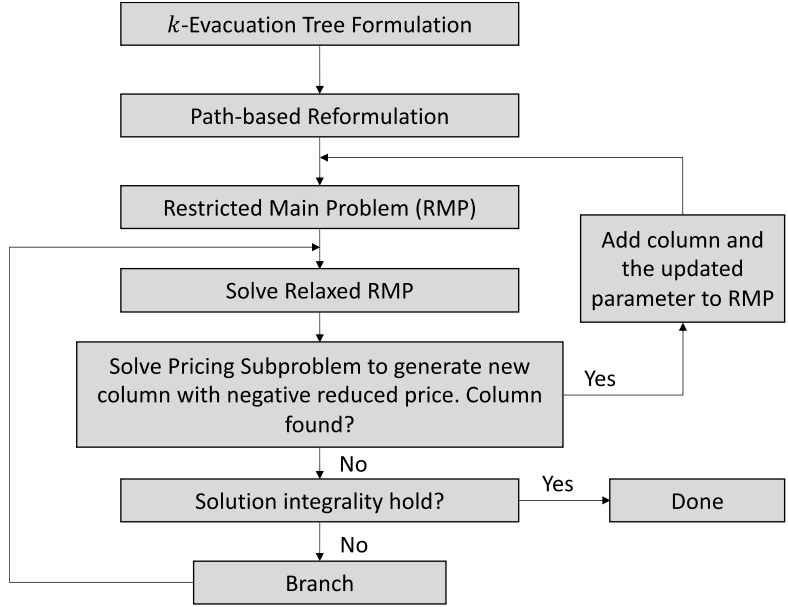

**Figure 4 Schema of the branch-and-price algorithm, adapted from** (Barnhart et al., 1998).

### Branching Scheme

The optimal solution of the Relaxed RMP may not satisfy the integrality constraints (14) and (23). In this case, we need to branch. There are several branching schemes that have been successful in the literature (Barnhart et al., 1998; Fischetti et al., 2002), including branching the most fractional variables or branching the most weighted variable. In the proposed solution method, we branch according to the following.

Let $\left(\tilde{x}^k, \tilde{\lambda}^{ko}\right)$ be the optimal fractional solution of the current relaxed RMP for all $o \in O, k \in K$. First, we branch the most fractional $\tilde{\lambda}^{ko}$ and create two new nodes: one corresponding to $\tilde{\lambda}^{ko}_{p*} = 0$ and another for $\tilde{\lambda}^{ko}_{p*} = 1$. Second, we adopt the branching strategy for the arcs from Tilk and Irnich (2018). We select one fractional $\tilde{x}^k_{ij}$ that maximizes $\sum_{(i,j) \in A} t_{ij} (\hat{v}_{ij}) x^k_{ij}$ and set $\tilde{x}^k_{ij*} = 1$ and tentative branching $\tilde{x}^k_{ij*} = 0$. Third, we guarantee that each branch node is not violating the contraflow constraints (7). If following the branching rule results in violating constraints (7), we drop that branch. These rules have been evaluated in past literature and have been shown to balance the branch-and-bound tree and (on average) to outperform other rules such as choosing the most fractional variables (Tilk and Irnich, 2018).

## NUMERICAL EXPERIMENTS

Experiments illustrate the application of the proposed evacuation network model on the Sioux Falls and South Florida transportation networks. Our analysis uses the Sioux Falls transportation network as the basis network to uncover the importance of parameters such as the driving range and the number and location of refueling stations deployment in the evacuation planning of alternative fuel vehicles. Then, we use the South Florida network to demonstrate the application of the $|K|$-evacuation tree routes planning problem in a larger network. The optimization problems (RMP and SP) and solution algorithms are modeled and solved using the commercial solver Gurobi (Gurobi Optimization, LLC, 2021) and its Python interface. The networks and their operations are modeled using NetworkX (Hagberg et al., 2008).

### Network Description

The evacuation problem is solved for the transportation network topology of the Sioux Falls and South Florida. The Sioux Falls network consists of 24 nodes and 76 links. We consider three types of network nodes, i.e., evacuee demand nodes, refueling station nodes, and a single shelter located at node 2 (the northeastern-most node). The alternative refueling stations are installed in nodes 5, 11, 12, 15, and 16. These positions were optimal deployment in the study of He et al. (2014), when solving the electric vehicle network equilibrium problem for habitual travel in the Sioux Falls network. In our analysis, we use these nodes as set of refueling stations and consider combinations from this set to generate different scenarios of refueling topology.

All transportation network properties (including the road free-flow travel times and links' capacities) are taken from data that is openly accessible (Transportation Networks for Research Core Team, 2021). **Figure 5** presents the transportation network topology of Sioux Falls, which is used in our numerical experiments. In a real-world evacuation, the evacuee demand originating from each node/zone could be retrieved from data collected via surveys, traffic detectors, GPS traces, or social media (e.g., Feng et al., 2020; Maas et al., 2019; Metaxa-Kakavouli et al., 2018; Sarma et al., 2020). However, in this paper, we make assumptions for the total evacuee demand from each evacuee node with total network flow of 360,600 vehicles to demonstrate the proposed evacuation problem, as shown in **Table 2**.

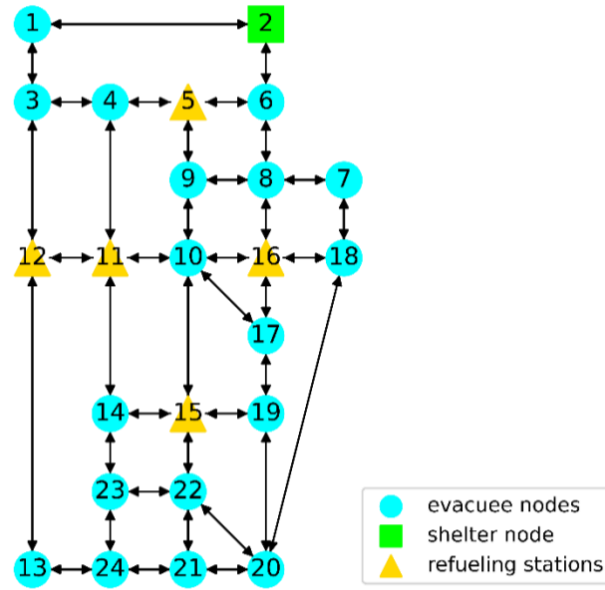

**Figure 5 Sioux Falls transportation network with refueling stations at nodes 5, 11, 12, 15, and 16.**

**Table 2 Evacuation Demand to Shelter Node 2**

| Node | Demand (vehicles) | Node | Demand (vehicles) | Node | Demand (vehicles) |
|------|------|------|------|------|------|
| 1 | 8800 | 9 | 16200 | 17 | 23400 |
| 2 | 0 | 10 | 45200 | 18 | 4800 |
| 3 | 2800 | 11 | 22300 | 19 | 12800 |
| 4 | 11600 | 12 | 13900 | 20 | 18500 |
| 5 | 6100 | 13 | 14600 | 21 | 11000 |
| 6 | 7600 | 14 | 14100 | 22 | 24400 |
| 7 | 12100 | 15 | 21400 | 23 | 14500 |
| 8 | 16700 | 16 | 26100 | 24 | 7700 |

The topology and properties of the South Florida transportation network are obtained from Sun et al. (2020). The South Florida network consists of 82 nodes and 234 links. The total travel demand is 65,707 vehicles with 83 OD pairs. We assume that the total evacuee demand of each evacuee node is the total flow originating from origin nodes of the habitual origin-destination (OD) matrix of the South Florida network. There are 31 nodes with zero evacuee travel demand, called "intersection nodes." Alternative refueling stations are installed in 30 nodes of the network. These nodes were the optimal deployment findings in the Sun et al. (2020) study, when solving the electric charging station problem for habitual travel in the South Florida network. In our analysis, we use these nodes as set refueling stations and consider combinations from this set to generate different scenarios of refueling topology. Additionally, we consider 5 shelters located at nodes 1, 50, 64, 41, and 81. We added node *s* as a dummy shelter node and 5 direct links to node *s* with zero travel time and infinite capacity. The South Florida numerical experiments are conducted using the modified transportation network where evacuees are all directed to the shelter node *s*. **Figure 6** presents the original South Florida network configuration (a) and the modified transportation network (b).

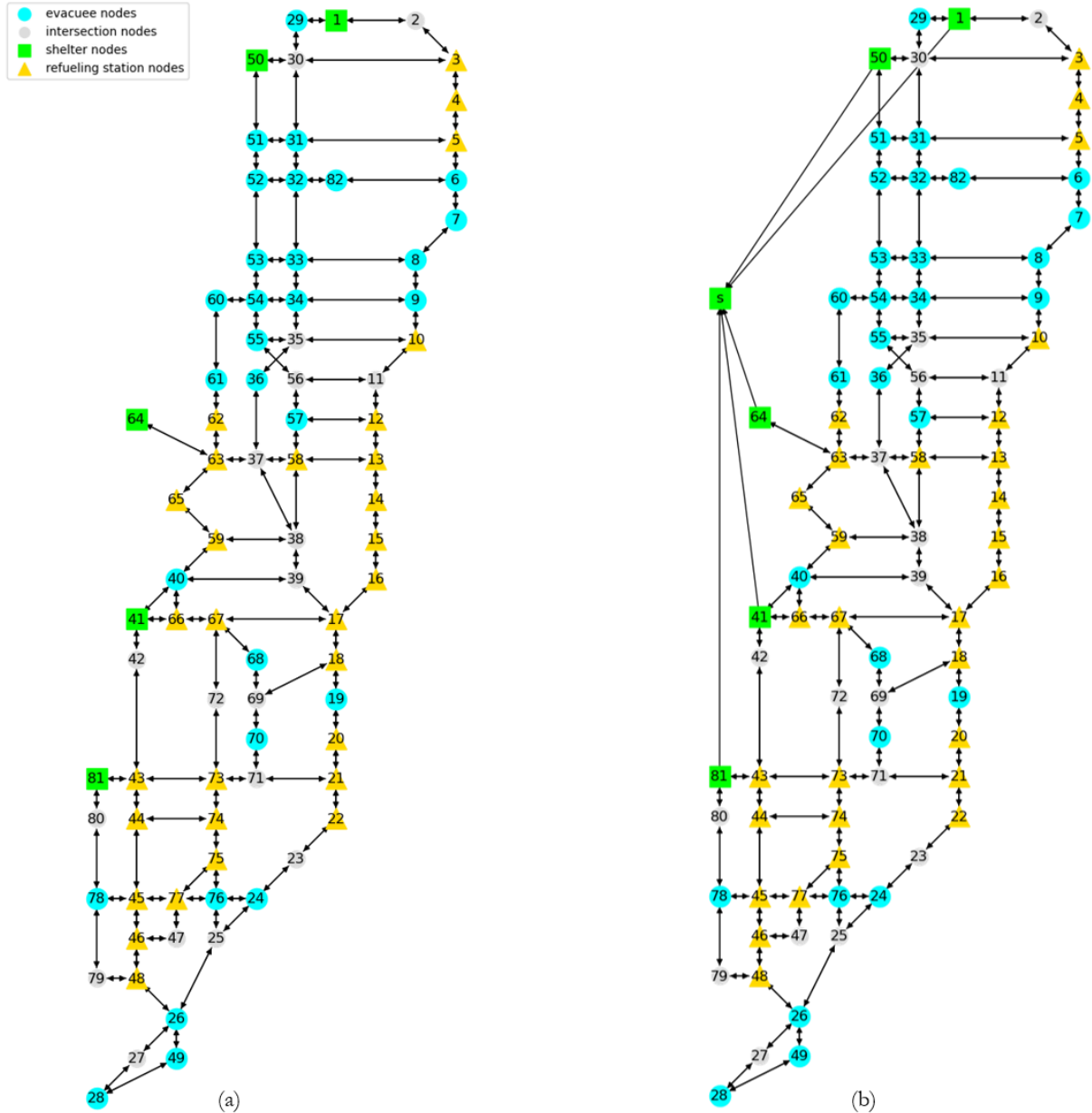

**Figure 6 (a) The South Florida transportation network with 30 refueling station nodes and 5 shelter nodes, (b) the modified South Florida transportation network with a safety node *s*.**

## Incorporating Refueling Needs in the Evacuation Route Plan

Initially, we examine the impact of incorporating the refueling needs of alternative fuel vehicles in the evacuation network planning and show the difference between the conventional evacuation route plans and the proposed evacuation route plan for alternative fuel vehicles. In this experiment, we solve the problem where every evacuee uses a single vehicle fuel type to solely highlight the refueling needs impact without the existence of multiple vehicle types, thus we set $|K| = 1$. Additionally, we solve the problem for two cases: (a) a conventional evacuation route plan, where vehicle range and refueling requirements are neglected, and (b) the proposed evacuation network problem where vehicles have a driving range, $\tau_k = 4$ hops. For the first case, the conventional vehicle evacuation tree problem adheres only to seamlessness and contraflow principles (Achrekar and Vogiatzis, 2018; Hasan and Van Hentenryck, 2020), which can be modeled using the subset of our formulation that corresponds to the objective function and constraints (1)—(7), (14), (17), and (18). For the second case, the proposed evacuation tree problem with refueling needs for alternative fuel vehicles is formulated as in (1)—(18). Then, we set the refueling stations' location as in **Figure 5** and the refueling time rate, $r_c^k = 15$ minutes per hop.

**Figure 7** shows the optimal evacuation route plans for conventional and alternative fuel vehicles. Both evacuation routes have a tree structure that meets the seamlessness evacuation attribute, but the resultant routes are different. In **Figure 7a**, the conventional evacuation plan guides every vehicle evacuee directly to the shelter without guaranteeing access to the refueling stations for every origin node. For example, every vehicle evacuee originating from nodes 18, 19, and 20 follows the 19-20-18-7-8-6-2 path (6 hops length) to the shelter directly without any access to refueling stations. However, when every evacuee uses a vehicle with a driving range proxy of $\tau_k = 4$ hops, the aforementioned evacuation path is not feasible since its length exceeds the driving range barrier and results in stranded and vulnerable evacuees from origin nodes 18, 19, and 20 during the evacuation process.

In contrast, our proposed evacuation plan, presented in **Figure 7b**, guarantees refueling access to each alternative fuel vehicle originating from nodes that need to detour and meet their vehicles' refueling requirements. The green highlighted area represents all evacuee origin demands that require access to refueling stations due to limited driving range. The proposed evacuation route plans indicate 12 evacuee origins with refueling needs during the evacuation process, which the conventional evacuation route plan does not capture. Focusing on evacuee demand from nodes 18, 19, and 20 in **Figure 7b**, our proposed evacuation route provides refueling access for each originating demand to meet their refueling needs prior to reaching the shelter, which could not be achieved in the conventional evacuation plan. We show that the conventional route is not feasible for alternative fuel vehicles with driving range barriers.

Furthermore, at the system level, we observe that the total evacuation time increases by 7.32% when incorporating the vehicles' refueling needs. This is due to the increased travel time caused by detours to the refueling stations and delays due to refueling. Our result indicates that the current conventional evacuation plan would undermines the driving range barriers and their impact on the evacuation performance, making access to shelter infeasible and leaving alternative fuel vehicle travelers vulnerable.

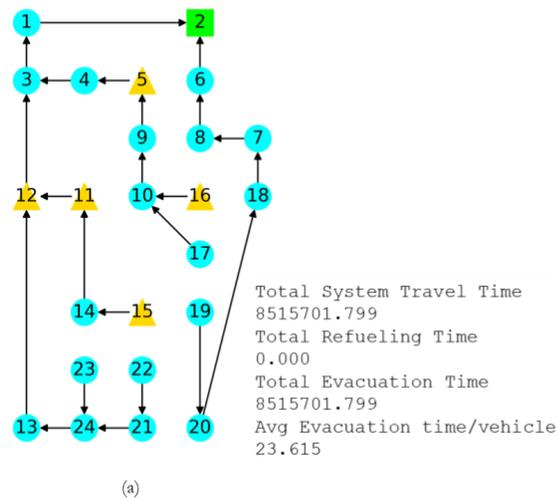

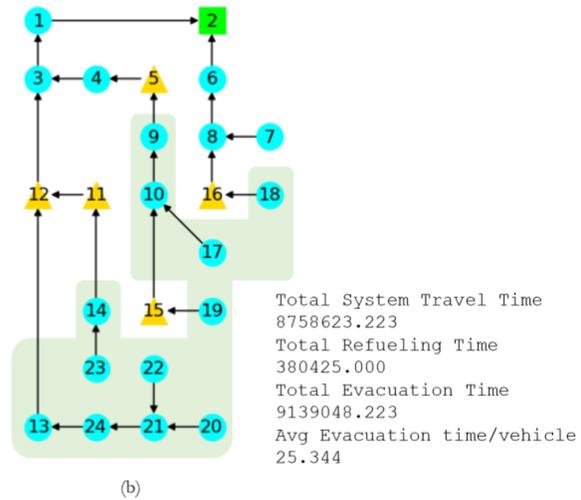

Note: time unit is in hours

(a) Total System Travel Time 8515701.799 / Total Refueling Time 0.000 / Total Evacuation Time 8515701.799 / Avg Evacuation time/vehicle 23.615

(b) Total System Travel Time 8758623.223 / Total Refueling Time 380425.000 / Total Evacuation Time 9139048.223 / Avg Evacuation time/vehicle 25.344

**Figure 7 Optimal evacuation tree routes for (a) conventional route plan; (b) proposed alternative fuel vehicles evacuation tree with refueling stations located at nodes 5, 11, 12, 15, 16 and driving range $\tau_k = 4$ hops.**
*Note: the green area indicates all evacuees' origins for which routes need to provide access to a refueling station*

**Driving Range Constraints Impact on the Evacuation Route Plan**

We further examine the impact of the vehicle driving range on the optimal evacuation route plan. We again consider $|K| = 1$ while the refueling station locations are as in **Figure 5** and have the same refueling efficiency, $r_c^k = 15$ minutes per hop.

**Figure 8** illustrates the optimal evacuation routes for four different vehicle types with driving range limits of $\tau_k = \{3, 5, 7, 8\}$ hops. The evacuation routes differ for different driving range bounds. The evacuation route of **Figure 7b** and **Figure 8a** show that as the vehicle's driving range is restricted to 3 hops, the total evacuation time increases significantly by 276.39%. There are 14 origin nodes whose evacuee demand needs to refuel before reaching the shelter. In **Figure 7b**, we see how evacuees originating from node 4 follow the 4-3-1-2 path without refueling needs, while the same evacuees need to refuel at the refueling station of node 5 following the 4-5-6-2 path in **Figure 8a**. Specifically, evacuees from nodes 7 and 18 are subject to major detours to meet their refueling needs at the station of node 12 instead of that on node 16. This is due to the system optimum (SO) objective of the evacuation traffic assignment. Based on the network properties, links (18,20) and (1,2) have greater road capacities than links (18,16) and (6,2) (Transportation Networks for Research Core Team, 2021); hence, the proposed evacuation route assigns evacuees on links (18,20) and (1,2) instead of links (18,16) and (6,2).

In **Figure 8b**, as expected, we observe improved evacuation performance with reduced numbers of origin nodes needing their evacuees to refuel, for the case of $\tau_k = 5$ hops. As the driving range of alternative fuel vehicles increases, the evacuation travel time decreases leading to better performance by 8.82%. Compared to the case of $\tau_k = 4$ hops, evacuees originating from node 18 do not require to reroute to the refueling station because the shortest evacuation path satisfies the current driving range barriers and becomes the optimal evacuation path. In addition, evacuees originating from node 13 do not require to reroute to the refueling station, however, we observe that this evacuee demand is passing through the station on node 12 to maintain the tree structure of the route.

Comparing **Figure 7a, Figure 8c, and Figure 8d**, the evacuation networks for $\tau_k = 7$ hops and $\tau_k = 8$ hops produce the same evacuation tree shape as the conventional plan. When $\tau_k = 7$ hops, the optimal route still requires evacuees from nodes 17 and 22 to refuel, even while following the shortest evacuation path to the shelter; while when $\tau_k = 8$ hops evacuees can directly reach the shelter due to their sufficient driving range. We show that the driving range constraint is not a significant barrier anymore in the evacuation route design and evacuees could traverse the shortest evacuation path. Our result indicates that as the vehicles' driving range increases, evacuees could fulfill their refueling needs while avoiding major detours.

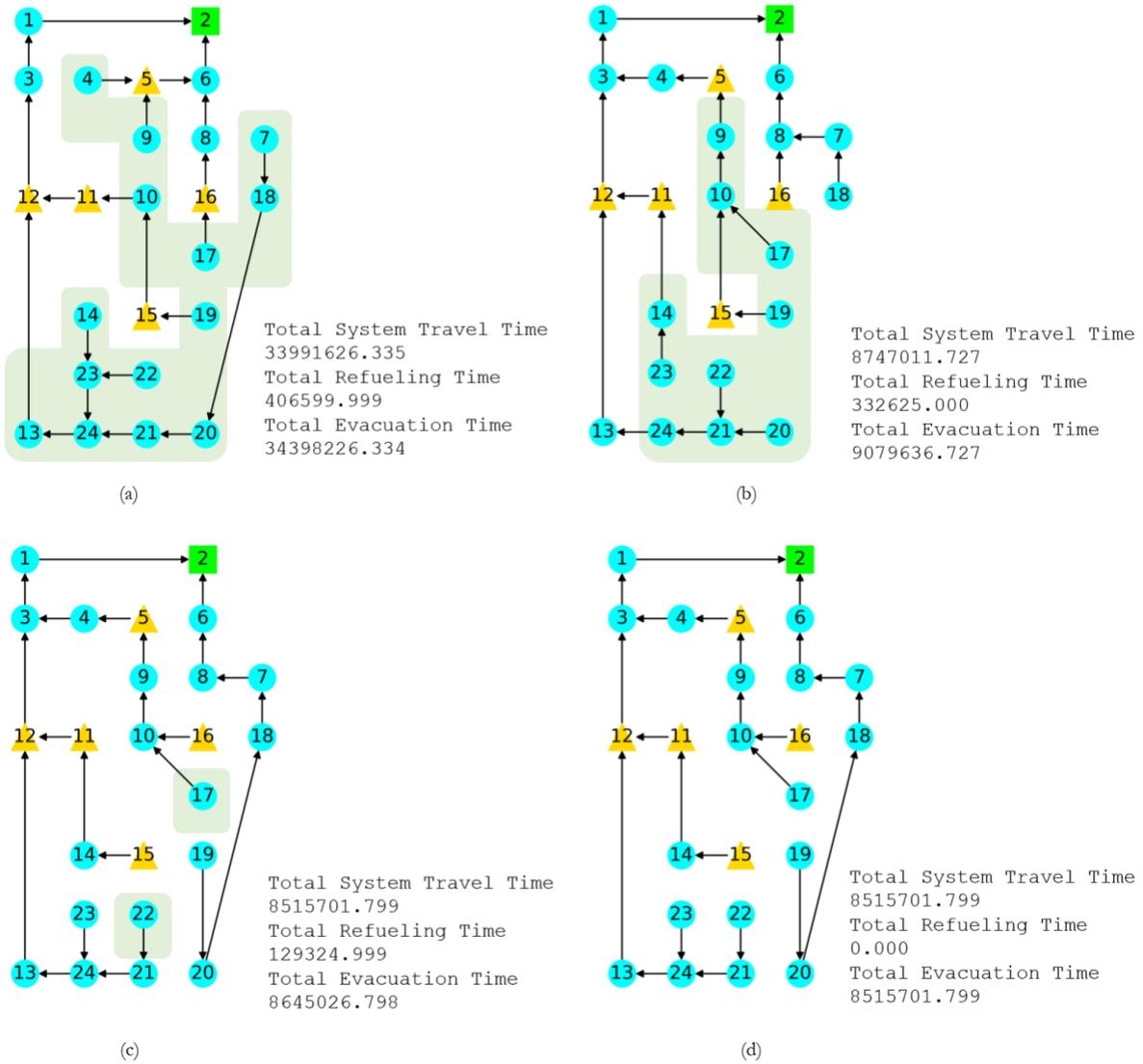

**Figure 8 Optimal evacuation tree routes with refueling stations located at nodes 5, 11, 12, 15, 16 for four cases of vehicles' driving range: (a) $\tau_k$ = 3 hops, (b) $\tau_k$ = 5 hops, (c) $\tau_k$ = 7 hops, and (d) $\tau_k$ = 8 hops.**
*Note: the green area indicates all evacuees' origins for which routes need to provide access to a refueling station*

In the next step, we evaluate the overall impact of driving range limit to the evacuation route design by enumerating extensive cases of driving ranges from $\tau_k = [0, 15]$. **Figure 9** presents evacuation performance results when the vehicle driving range varies from 0 to 15 hops. As the alternative fuel vehicle driving range increases, both the travel and refueling times improve, as expected. We confirm that vehicle range impacts traffic assignment performance, a topic also studied in prior literature (Erdoğan and Miller-Hooks, 2012; He et al., 2014; Jiang et al., 2012). Our numerical experiments show that the range of alternative fuel vehicles also impacts the evacuation performance.

For $\tau_k \geq 8$ hops, the total evacuation time of an alternative fuel vehicle type converges to the evacuation time of the conventional vehicles' plan. This shows that having a large driving range would enable evacuees to traverse the shortest evacuation path without any need of refueling detours. The conventional evacuation plan is not feasible for an alternative fuel vehicle with range $\tau_k < 8$ hops. This shows the necessity for emergency coordinators to design unique evacuation routes for alternative fuel vehicles due to their low driving range compared to conventional gasoline vehicles.

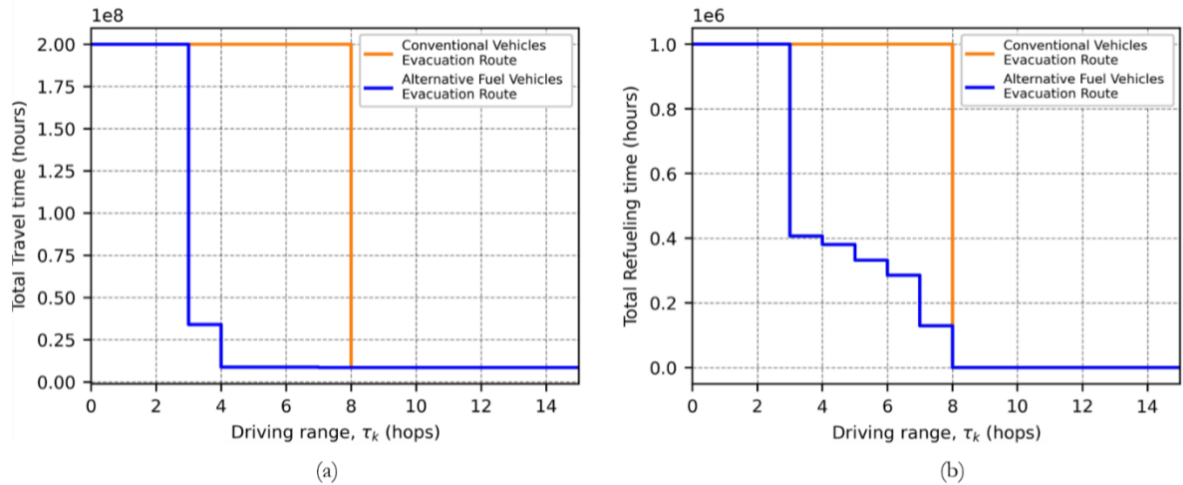

**Figure 9 (a)** The evacuation travel time and **(b)** the recharging time of an alternative fuel vehicle evacuation in the Sioux Falls network with charging stations located at nodes 5, 11, 12, 15, and 16 and the driving range varying from 0 to 15 hops.

Furthermore, the case of $\tau_k = 0$ is infeasible since vehicles do not have the ability to traverse any of the network links due to their extremely short driving range. When $\tau_k > 0$, infeasibility may still occur due to the driving range bound, the sparse refueling infrastructure network, and the route's desired seamlessness attribute. For example, for the overly restrictive driving range of $\tau_k = 1$ hop, evacuees need to detour for refueling even one hop away from the shelter node. Hence, if there is no refueling stations located on nodes 1 and 6 (the two nodes adjacent to the shelter), then no feasible tree rooted at the shelter exists. We consider these infeasible cases for validation purposes of the evacuation model.

**Refueling Network Density Impact on Evacuation Route Plan**
This section explores the impact of refueling station topology on the evacuation routing plan. We consider a single vehicle type ($|K| = 1$) with the same refueling rate of the vehicle, $r_c^k = 15$ minutes per hop. For the refueling stations deployment, we generate combinations from the set of open refueling stations and evaluate the cases of only having 1 station, 2, 3, 4, or all 5 stations available. We also vary the driving range from 0 to 15 hops by enumerating the scenarios for every possible $\tau_k$ value and every possible subset of those five stations' location. Hence, we review the results of $16 \cdot (2^5)$=512 scenarios.

**Figure 10** illustrates the designed evacuation routes for four indicative and different cases of refueling station deployment with driving range of vehicles as $\tau_k = 4$ hops. We observe that both the total number of installed refueling stations and their deployment position affect the evacuation routing plan and result in unique routing paths. For instance, comparing **Figure 10a and Figure 10b**, it can be observed that the evacuation routing plan and the overall performance are unique due to different refueling station siting, even with same total number of refueling stations on the network. In **Figure 10b,** we observe that as we site the refueling stations closer to the shelter, the total evacuation time improves significantly. This result is expected as we design routes to direct evacuees to the shelter node 2. The refueling station placement near the shelter would provide opportunities for evacuees to fulfil their required refueling needs using the shorter route without major detour.

Furthermore, the evacuation performance improves as we increase the number of installed refueling stations. When an additional refueling station is sited in node 16, the evacuation time improves by 39.34% when comparing cases in **Figure 10a and Figure 10c**. We observe 3.69 % improvement of the evacuation time by adding a refueling station at node 16, as shown in **Figure 10b and Figure 10d**. However, the operational improvement, as the number of installed refueling station increases, is not linear, and it depends on the position of the infrastructure. Comparing **Figure 10b and Figure 10c**, we show that a denser refueling station network would not solely guarantee a faster evacuation process. A strategic refueling siting plan, enabling connectivity with the shelter node, plays a significant role in prompt evacuations.

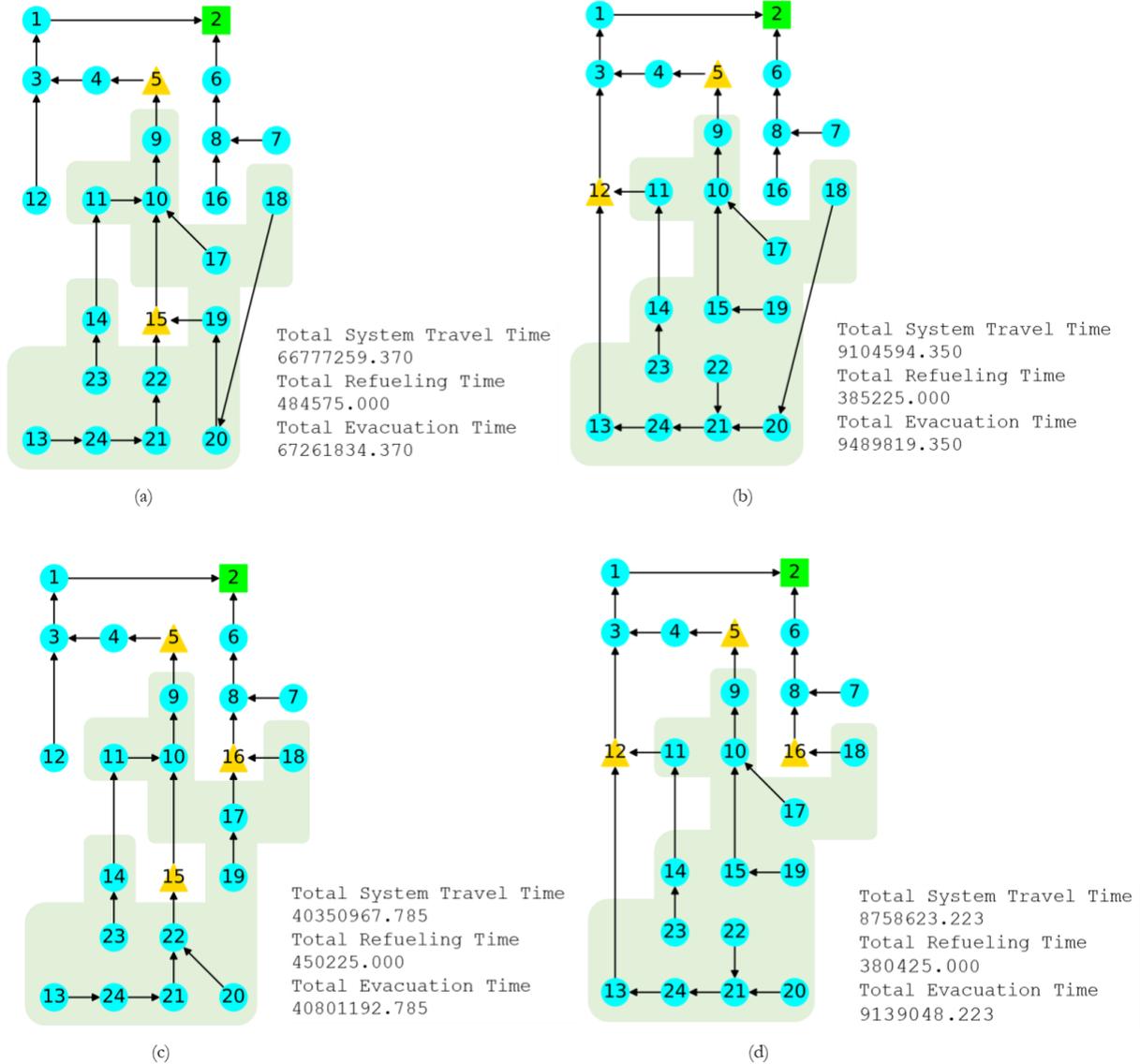

**Figure 10** Optimal evacuation tree routes with vehicle driving range $\tau_k$ = 4 hops, for four cases of refueling station deployment.
*Note: the green area indicates all evacuees' origins for which routes need to provide access to a refueling station*

To evaluate the overall impact of refueling stations density on the network, we generate combinations of the set of open refueling stations {5, 11, 12, 15, 16}, and evaluate the cases of only having 1 refueling station, 2, 3, 4, or all 5 refueling stations available. **Figure 11** shows the total evacuation time for all generated combinations of refueling station locations. Each sub-graph represents the evacuation performance given a different number of installed refueling stations. Each line in each subgraph represents the objective outcome for different values of $\tau_k$, given a refueling deployment scenario. The result confirms that different siting of refueling station produces unique evacuation operational performance, as expected.

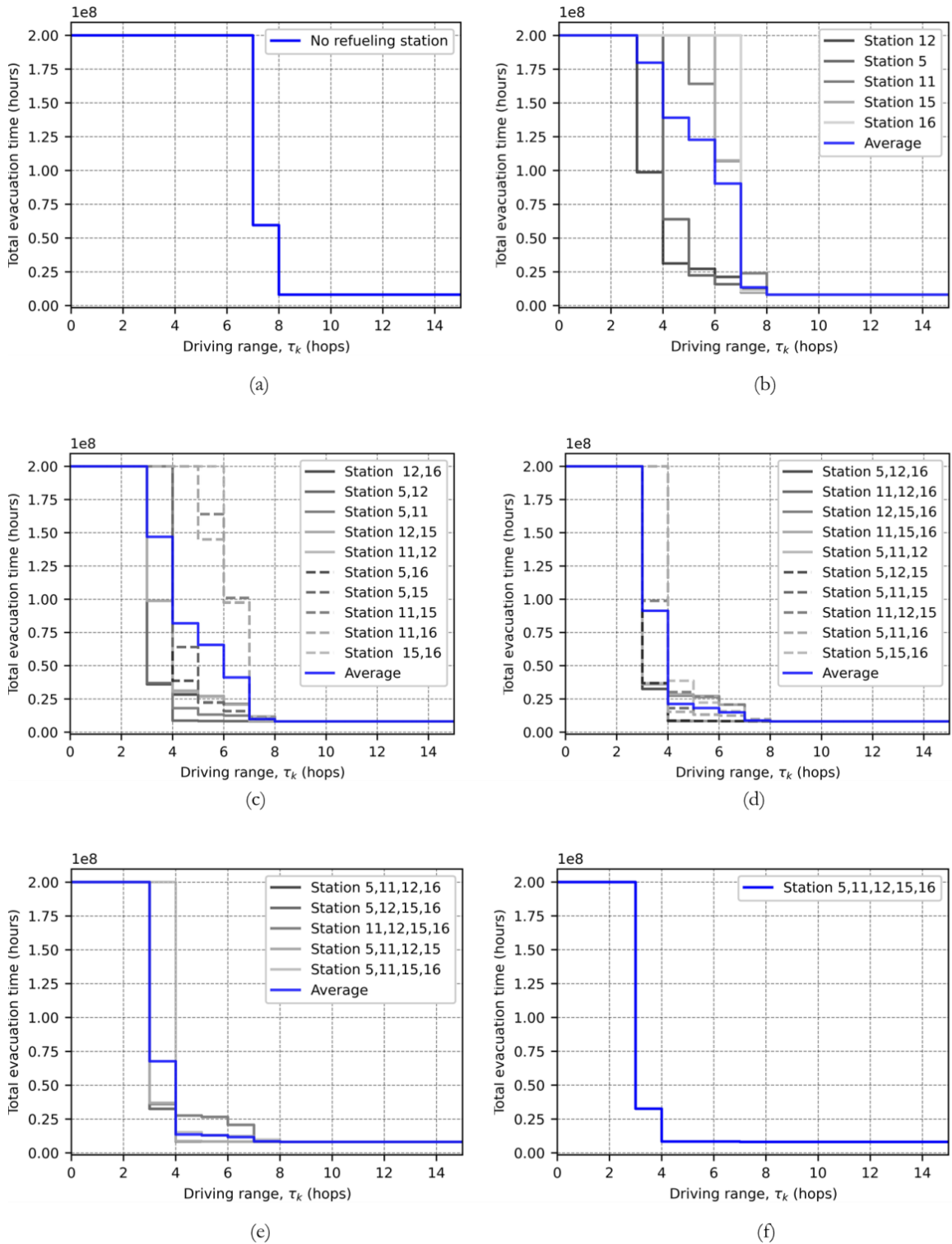

**Figure 11** Results of refueling station location analysis for the Sioux Falls evacuation network: scenario (a) corresponds to no charging infrastructure, scenario (b) includes every case of only 1 refueling station sited, scenario (c) every case of 2 refueling stations sited, scenario (d) every case of 3 refueling stations sited, scenario (e) every case of 4 refueling stations sited, and scenario (f) cases with 5 refueling stations sited.

**Figure 12** summarizes the impact of the refueling station density on the average evacuation routing plan performance for each refueling station density in the network. We observe that as the refueling station network becomes denser, the evacuation performance also, on average, improves. The case of 5 refueling stations provides consistently the best evacuation performance for different values of $\tau_k$. Focusing for $\tau_k \geq 8$ hops, the total evacuation time of alternative fuel vehicles again converges to the total evacuation time of a conventional evacuation plan, as shown in **Figure 7a**. This result indicates that the refueling station density would significantly affect the evacuation performance when vehicles have short driving ranges since, on average, denser infrastructure placement improves both evacuation and refueling times. When the vehicles' range is large, evacuees could traverse the shortest evacuation path without need of detour.

This result confirms a similar finding noted in literature (Erdoğan and Miller-Hooks, 2012) that has examined the impact of recharging network density and stations locations on the routing assignment and its feasibility. We verify that the dependency on refueling infrastructure also holds in the evacuation context. We conclude that a denser alternative refueling network improves the average evacuation performance, but the siting of stations matters and provides unique evacuation performance. This poses new challenges for emergency planners regarding how to strategically coordinate with industry for alternative refueling station deployment that supports evacuations, when such siting decisions are often not made by public agencies.

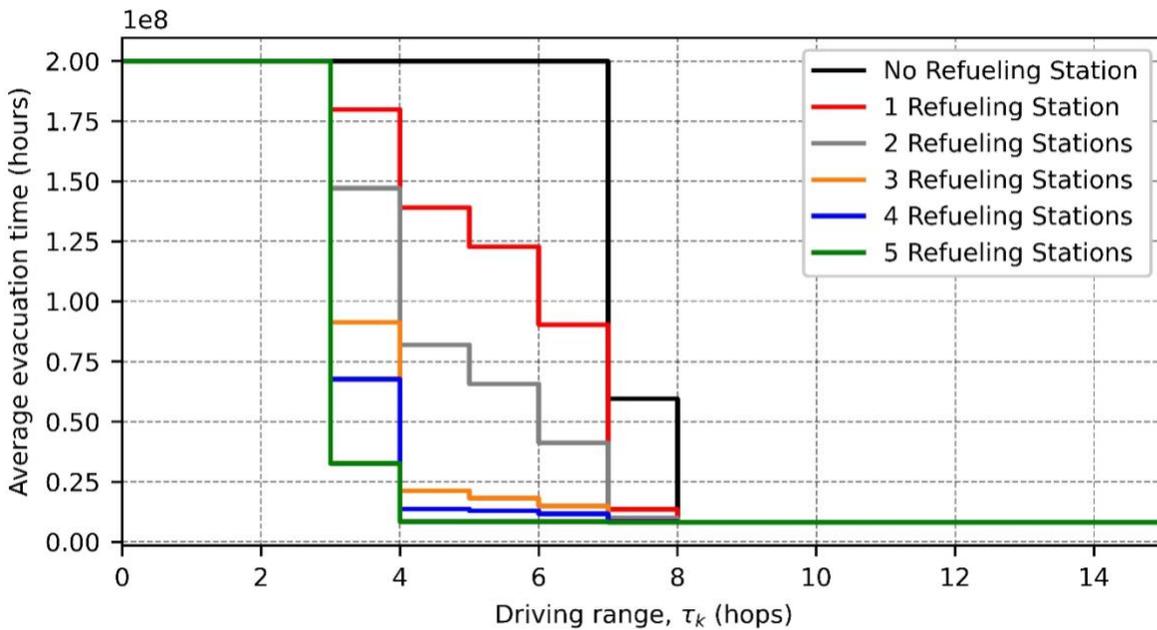

**Figure 12 Average evacuation performance that corresponds to scenarios of different numbers of operating refueling stations on the evacuation network.**

## $|K|$-Evacuation Route Plans

### *Application in the Sioux Falls Transportation Network*

Now, we examine the impact of coordinating multiple alternative fuel vehicles' evacuation planning. We assume three types of alternative fuel vehicles: vehicle type I represents an automobile with a large driving range and a dense refueling infrastructure network, vehicle type II represents an automobile with a short driving range and a sparse refueling network, and vehicle type III represents a vehicle fuel type with a medium driving range and a sparse refueling network. The vehicle characteristics and refueling station deployment for each type are presented in **Table 3**. We examine three cases of evacuation route design: a single vehicle fuel type, two vehicle fuel types, and three vehicle fuel types. In case A, we consider that every evacuee uses the same vehicle type with a large driving range and a dense refueling network. Note that the evacuation route for case A has been presented in **Figure 8d**.

**Table 3 The Vehicle Characteristics and Refueling Station Position for Cases A, B, and C**

| Case | Number of Vehicles, $|K|$ | Vehicle Type | Population Share | Refueling Station Position | Driving Range Limit, $\tau_k$ | Refueling Rate, $r_c^k$ |
|------|------|------|------|------|------|------|
| A | 1 | I | 100% | 5,11,12,15,16 | 8 | 5 |
| B | 2 | I | 50% | 5,11,12,15,16 | 8 | 5 |
|  |  | II | 50% | 5,15,16 | 4 | 30 |
| C | 3 | I | 40% | 5,11,12,15,16 | 8 | 5 |
|  |  | II | 30% | 5, 15, 16 | 4 | 30 |
|  |  | III | 30% | 5, 11, 16 | 6 | 10 |

In case B, we consider two alternative fuel vehicle types, (i) a vehicle fuel type with a large driving range, dense refueling stations network, and short refueling time, and (ii) a vehicle fuel type with short driving range, a sparse refueling network, and long refueling time. We also assume that the vehicle population is shared equally to the evacuee demand in the network. **Figure 13** presents the evacuation tree plan for each vehicle type in case B. We observe that each evacuation route successfully maintains the principle of seamlessness and contraflow to avoid conflicts with other vehicle fuel types, while fulfilling their refueling needs. Comparing **Figure 8d** and **Figure 13a**, the evacuation route of evacuees with vehicle type I changes to avoid the traffic conflict with evacuees of vehicle type II. For instance, in case A (see **Figure 8d**), the evacuee demand originating from node 21 is routed to node 24 and evacuees from node 24 are routed to node 13. This traffic assignment is conflicting with the routing assignment of evacuees of vehicle fuel type II from node 13 and 24, as shown in **Figure 13b**. Our results produce alternative evacuation plans for vehicles of type I to avoid the conflicting assignment between those demand nodes, as in **Figure 13a**.

The total evacuation time for case B increases significantly by 142.50%, compared to case A. The evacuation performance declines due to the consideration of alternative fuel vehicle barriers of vehicle type II. When the driving range and refueling station density increases, the evacuation performance improves. In case B, vehicles of type II have short driving ranges, a sparse refueling infrastructure network, and long refueling times, leading to significantly reduced evacuation operation performance compared to case A.

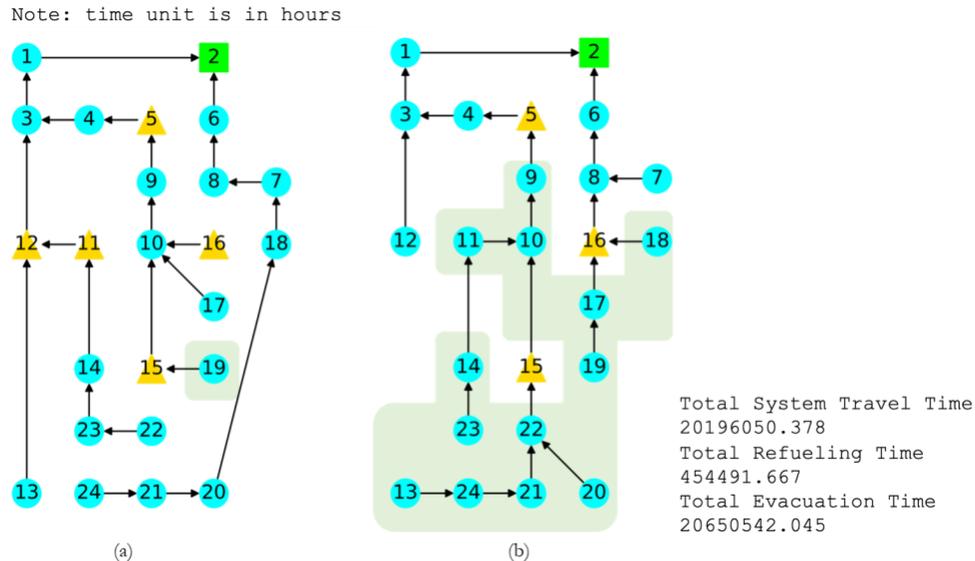

**Figure 13 $|K|$-evacuation route plan for Case B with 2 alternative fuel vehicle types: (a) vehicle type I, and (b) vehicle type II.**

We also consider three alternative fuel vehicle types: we introduce a vehicle fuel type III with moderate driving range and a sparse refueling station network. We assume that the vehicle type population of the evacuee demand is 40%, 30%, and 30%, respectively. **Figure 14** presents the evacuation tree plan for each vehicle type in case C. Comparing **Figure 8d**, **Figure 13a**, and **Figure 14a**, the resultant evacuation route for vehicle type I is changing again to accommodate the refueling needs of both vehicle types II and III. The evacuation route plan for vehicle type II is also changing, as we consider vehicle type III, compared to **Figure 13b** and **Figure 14b**. We observe that the total evacuation time of case C increases by 101.84% compared to the case A, but improves by 20.14% compared to case B. The result is expected since vehicle type III has a greater range compared to vehicle type II and given the vehicle population distribution on the network.

Total System Travel Time
16853616.982
Total Refueling Time
334821.667
Total Evacuation Time
17188438.649

**Figure 14** $|K|$-evacuation route plan for Case C with 3 alternative fuel vehicle types: (a) vehicle type I, (b) vehicle type II, and (c) vehicle type III.

As we consider multiple vehicle fuel types, the evacuation route design becomes more complex. Accounting for each vehicle characteristics and refueling stations availability while maintaining the attributes of seamlessness and contraflow leads to a slower overall evacuation performance. At the same time, it ensures that different fuel vehicle types can reliably access their refueling stations. We confirm that the previous insights of driving range impacts and refueling station density when $|K| = 1$ also apply to the multiple alternative fuel vehicle evacuation route design. We also observe that the evacuation route for each vehicle fuel type could be unique as we consider the heterogeneity of vehicle characteristics in the network. For example, the result shows that the evacuation route plan for vehicle type I are different for each case due to the heterogeneous share of vehicle fuel types in the transportation network. Considering the rise of various alternative fuel vehicle types' adoption in the near future (Davis and Boundy, 2021), this finding points out the pressing need for emergency planners to strategically manage and guarantee an efficient and safe evacuation routing plan for diverse vehicle fuel types.

### Application in the South Florida Transportation Network

In this section, we apply the $|K|$-evacuation tree route plan problem to the South Florida network. We again assume three types of alternative fuel vehicles: vehicle type I represents an automobile with a large driving range and a dense refueling infrastructure network, vehicle type II represents an automobile with a short driving range and a sparse refueling network, and vehicle type III represents a vehicle fuel type with a medium driving range and a sparse refueling network. The vehicle characteristics and refueling station deployment for each type are presented in **Table 4**. We examine the same three cases of evacuation route design: a single vehicle fuel type, two vehicle fuel types, and three vehicle fuel types. **Table 5** describes the refueling station deployment for each vehicle fuel type.

**Table 4 The Vehicle Characteristics and Refueling Station Position for Cases D, E, and F**

| Case | Number of Vehicles, $|K|$ | Vehicle Type | Population Share | Number of Refueling Station | Driving Range Limit, $\tau_k$ | Refueling Rate, $r_c^k$ |
|---|---|---|---|---|---|---|
| D | 1 | I | 100% | 30 | 6 | 15 |
| E | 2 | I | 50% | 30 | 6 | 15 |
|   |   | II | 50% | 15 | 4 | 30 |
| F | 3 | I | 40% | 30 | 6 | 15 |
|   |   | II | 30% | 15 | 4 | 30 |
|   |   | III | 30% | 22 | 5 | 15 |

**Table 5 The Refueling Station Deployment for Vehicle Type I, II, and III**

| Vehicle Type | Refueling Station Deployment Variation | Number of Refueling Stations | Refueling Station Position |
|---|---|---|---|
| I | Base (Sun et al., 2020) | 30 | 3,4,5,10,12,13,14,15,16,17,18,20,21,22,43,44, 45,46,48,58,59,62,63,65,66,67,73,74,75,77 |
| II | 50% random selection | 15 | 14, 73, 66, 5, 13, 17, 62, 45, 46, 22, 16, 12, 4, 77, 65 |
| III | 75% random selection | 22 | 22, 66, 10, 73, 77, 3, 62, 75, 65, 4, 18, 17, 12, 14, 59, 45, 48, 58, 74, 63, 67, 16 |

**Figure 15** presents the optimum evacuation route for case D where every evacuee uses the same vehicle type with a large driving range and a dense refueling network. We show that the tree structure is maintained rooted at the node $s$ with a driving range bound of 6 hops. We also observe that there are 30 evacuee nodes whose demand exits through shelter node 64, 26 evacuee nodes whose demand exits through shelter node 41, 9 evacuee nodes whose demand exits through shelter node 50, 7 evacuee nodes whose demand exits through shelter node 1, and 4 evacuee nodes whose demand exits through shelter node 81. Using the modified network, we show that our model can provide optimal shelter assignment for every evacuee vehicle that minimizes the total system's evacuation time.

**Figure 16** presents the evacuation tree plan for each vehicle type in case E. We consider two alternative fuel vehicle types, i.e., vehicle type I and vehicle type II, and the vehicle population is distributed equally in the network. We observe that each evacuation route successfully maintains the principle of seamlessness and contraflow to avoid conflicts with other vehicle fuel types, while fulfilling their refueling needs. Comparing **Figure 15** and **Figure 16a**, we consistently find that the evacuation route of evacuees with vehicle type I changes to avoid the traffic conflict with evacuees of vehicle type II, maintaining the contraflow and seamlessness principles. The total evacuation time for case E increases by 84.80%, compared to case D, due to the consideration of alternative fuel vehicle barriers of vehicle type II.

**Figure 17** presents the evacuation tree plan for each vehicle type in case F with three alternative fuel vehicle types. Comparing **Figure 15**, **Figure 16a**, and **Figure 17a**, the resultant evacuation route for vehicle type I is changing again to accommodate the refueling needs of both vehicle types II and III. Similarly, the evacuation route plan for vehicle type II is also changing, as we consider the addition of vehicle type III in the fleet composition, as shown in **Figure 16b** and **Figure 17b**. We observe that the total evacuation time of case F increases by 43.32% compared to the case D, and is improved by 22.48% compared to case E. In this section, we demonstrate the application of the $|K|$-evacuation tree route plan problem to a larger network. We show that the previous insights regarding the impact of vehicle driving range, refueling station density, and vehicle heterogeneity found in the Sioux Falls case, also apply in a larger transportation network.

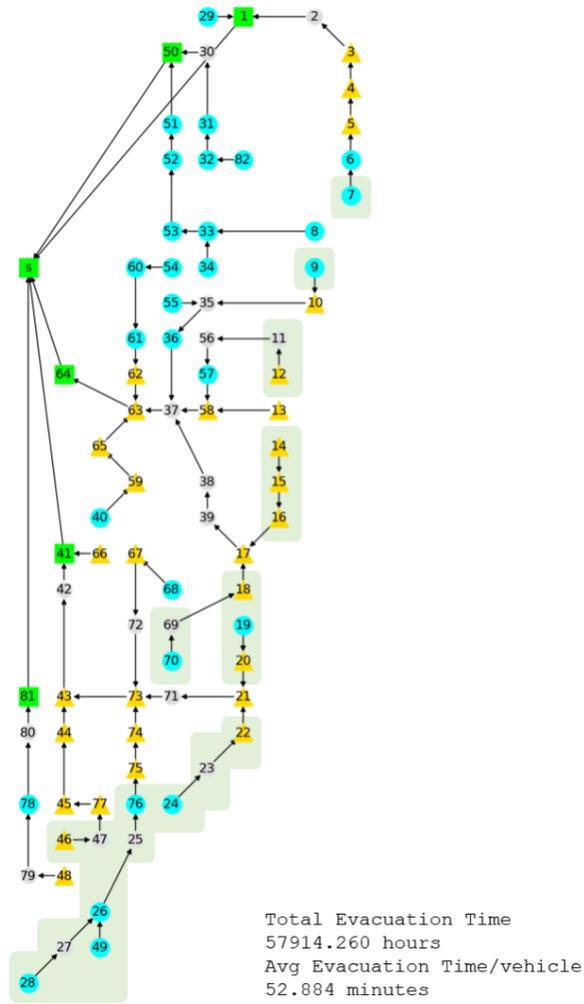

Total Evacuation Time
57914.260 hours
Avg Evacuation Time/vehicle
52.884 minutes

**Figure 15 Evacuation route plan for Case D with 1 alternative fuel vehicle type, vehicle type I**

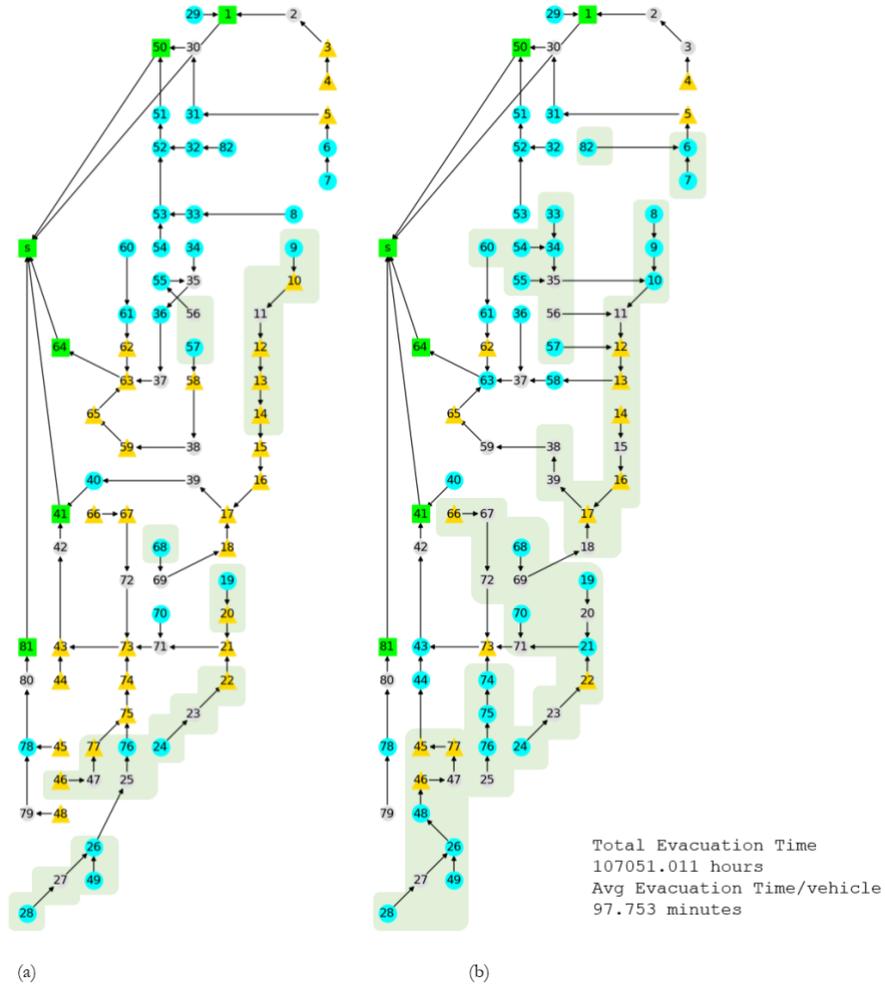

(a)                    (b)

Total Evacuation Time
107051.011 hours
Avg Evacuation Time/vehicle
97.753 minutes

**Figure 16 |$K$|-evacuation route plan for Case E with 2 alternative fuel vehicle types: (a) vehicle type I, and (b) vehicle type II.**

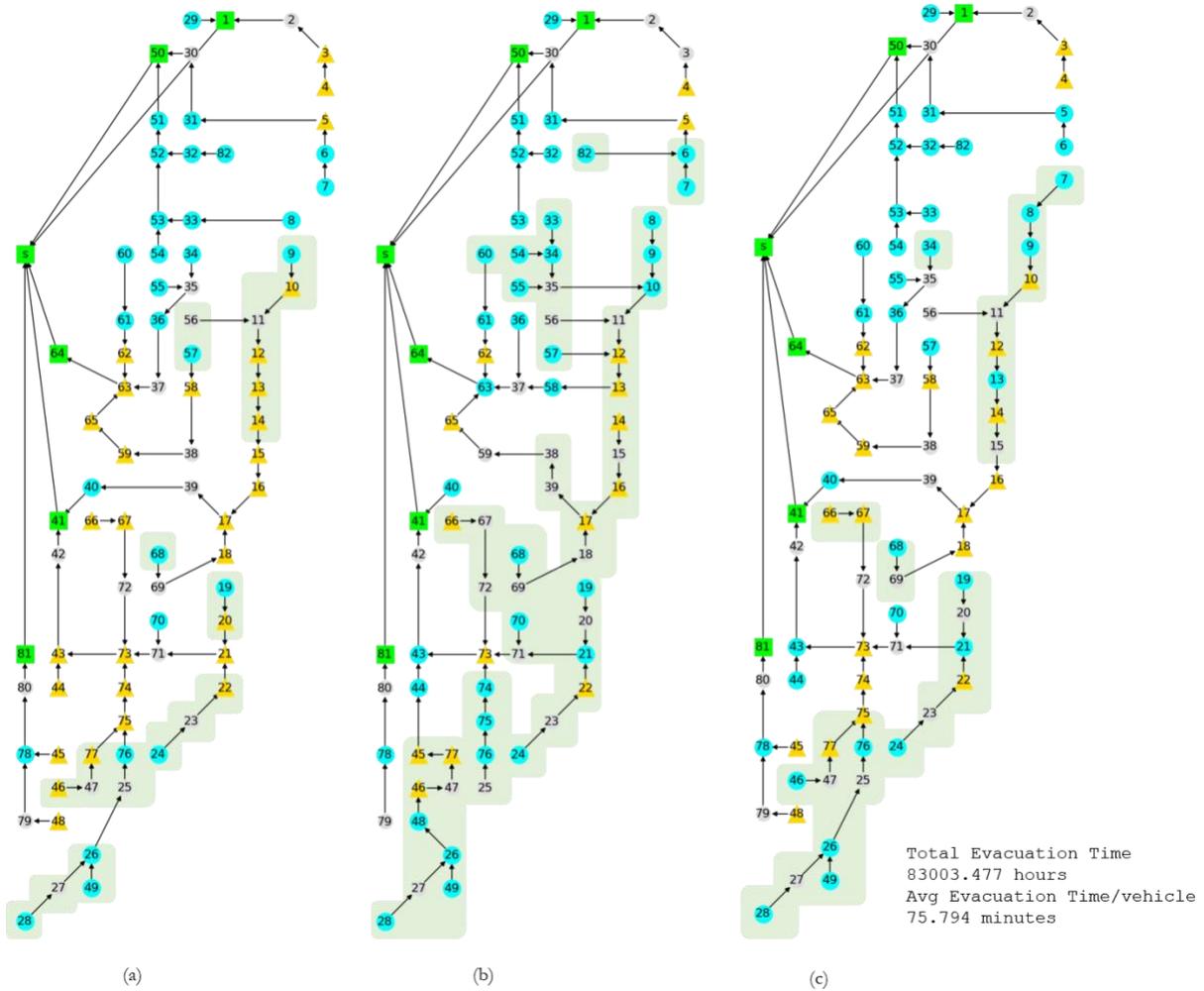

**Figure 17** $|K|$-evacuation route plan for Case F with 3 alternative fuel vehicle types: (a) vehicle type I, (b) vehicle type II, and (c) vehicle type III.

## Computational Performance

In this section we summarize the computational performance of our proposed evacuation problem and the matheuristic branch-and-price/column-generation. We compare the efficiency of the proposed solution method with the GUROBI solver result. The time limit is set to be 1500000 seconds. We run the model using a laptop with processor 1.4 GHz Quad-Core Intel Core i5 and memory 16 GB 2133 MHz LPDDR3.

**Table** 6 shows the summary of the computational performance of our numerical experiments. It shows that the branch-and-price/column-generation is effective in reducing computational time when solving the proposed evacuation tree problem for a set of alternative fuel vehicle types, $K$. For $|K| = 2$ vehicle fuel types (see **Figure 13**), the proposed approach solves the problem significantly faster, as opposed to Gurobi. Similarly, in the case of $|K| = 3$ vehicle fuel types (see **Figure 14**), our proposed solution method could result in determining the evacuation route in reasonable time, while the solver could not solve it (reached the time limit).

In terms of optimality, our matheuristic approach offers no guarantees on optimality gaps or on approximation. That said, we observe that in most of our instances the matheuristic returned a solution of high quality. More importantly in the evacuation planning setting, the matheuristic outperforms the commercial solver in all instances as far as computational time is concerned.

**Table 6 Summary of Solution's Computational Performance**

| Case Figure | Direct Solver GUROBI | | | Matheuristic Branch-and-Price/Column-Generation | | | | |
|---|---|---|---|---|---|---|---|---|
| | Objective (hour) | Avg Evacuation Time (hour) | Time (second) | Objective (hour) | Avg Evacuation Time (hour) | #Nodes | #Columns | Time (second) |
| 7b | 9139048.22 | 25.34 | 655 | 9246888.99 | 25.64 | 3 | 91 | 61 |
| 8a | 34398226.33 | 95.39 | 3041 | 34847198.82 | 96.64 | 6 | 288 | 278 |
| 8b | 9079636.73 | 25.18 | 1262 | 9270309.15 | 25.71 | 5 | 235 | 126 |
| 8c | 8645026.8 | 23.97 | 1828 | 8826572.36 | 24.48 | 6 | 385 | 216 |
| 8d | 8515702.8 | 23.62 | 2868 | 8600859.83 | 23.85 | 6 | 385 | 216 |
| 10a | 67261834.38 | 186.53 | 119 | 68338023.7 | 189.51 | 4 | 214 | 34 |
| 10b | 9489819.35 | 26.32 | 180 | 9660636.11 | 26.79 | 3 | 87 | 25 |
| 10c | 40801192.79 | 113.15 | 138 | 41821222.6 | 115.98 | 3 | 84 | 34 |
| 10d | 9139048.22 | 25.34 | 523 | 9239577.75 | 25.62 | 4 | 140 | 25 |
| 13 | 15794950.57 | 43.80 | 1224480 | 20650642.05 | 57.27 | 10 | 558 | 705 |
| 14 | 97671734.27 | 270.86 | 1500000* | 17188438.65 | 47.67 | 23 | 934 | 1286 |

*Terminated due to pre-determined time limit*

## POLICY IMPLICATIONS

In this section, we discuss the policymaking implications of the numerical experiments' findings and the implementation of our proposed evacuation model by emergency planning agencies. We provide recommendations and insights by highlighting the advantages of implementing the hop constraints approach for the evacuation planning of alternative fuel vehicles. Then, we discuss vehicle range and refueling density tradeoffs regarding evacuation route planning and point out how planners could strategically address these considerations. Finally, we address the policy implication of the diversity of vehicle fuel types for evacuation planning.

### Hop Constraints Opportunities in Evacuation Route Planning

In the proposed mathematical model, we implement hop constraints that approximate the distance or length of a route in the network with refueling limitations. Our hop constraints provide generality and efficiency in controlling the flow and the service level of the network at the macroscopic level. This approach enables us to use generic parameters to plan for a heterogeneous fleet of vehicles.

It would be challenging for emergency planners to incorporate the actual distances in the network as well as the driving range of each vehicle of fuel type $k$ in a preemptive evacuation planning. The hop constraints allow us to integrate the topology's distance and the variability of vehicle driving range in the same measurement unit. Instead of using the actual distance in the network and vehicle range variability among drivers, the evacuation route plan decision could be simplified to consider the controlled evacuation range as hop constraints in the network. The distance of one hop could represent the approximate controlled evacuation distance in the network used for evacuation planning. The hop constraint, $\tau_k$, could represent the driving range as a proxy of alternative fuel vehicles during evacuation planning; formulated as $\tau_k = \left\lfloor \frac{driving\ range}{controlled\ spacing} \right\rfloor$. We choose the floor function to capture the driver's comfortable range, which is usually represented as a certain percentage of the vehicle fuel storage volume or the preferred comfortable range. We assume that the link travel time of each hop would be divided proportionally by the total hops of the actual distance in the network.

**Figure 18** shows an example of converting actual distance into a controlled evacuation hop range. Consider a toy network with the information of link length in a miles unit and link travel time in a minute unit. We assume that the controlled evacuation spacing is within 100 mi; hence, one hop would represent any distance within 100 mi in the network. For instance, the distance between node 1 and the shelter, 325 mi, is approximated as 4 hops. We then create dummy nodes 3, 4, and 5 with zero demand to maintain hop distance in the network. Similarly, the distance between the refueling station and the shelter (200 mi) is approximated as 2 hops with the addition of dummy node 6 with zero demand. **Figure 18b** illustrates the modified hop-distance network. Consider a Tesla X vehicle with 333 mi of range. In the evacuation model, this driving range is equal to 3 hops and the feasible paths for this vehicle are 1-2-s and 1-r-6-s, where the optimal path depends on the minimized travel time and refueling time. For a Chevrolet Bolt EV with 259 mi of range, its range is equal to 2 hops. The feasible evacuation route plan is path 1-r-6-s.

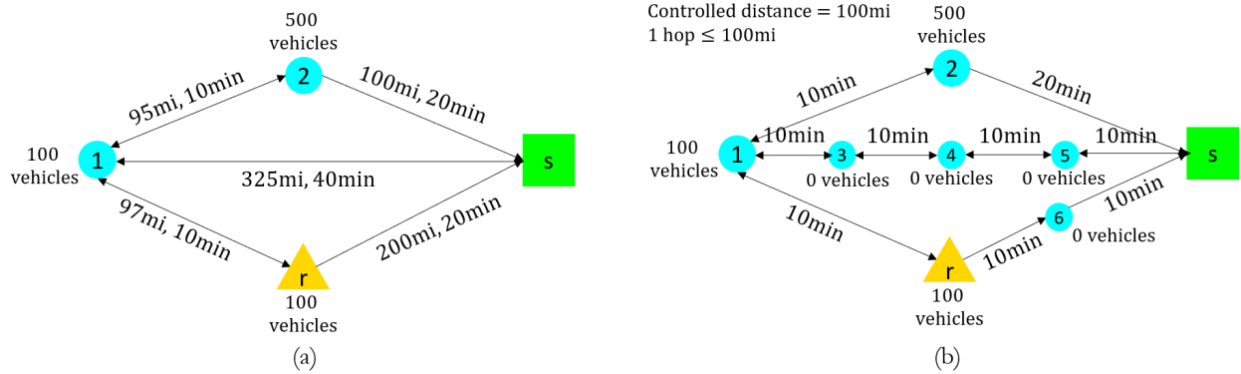

**Figure 18 Example of (a) network topology with link distance parameters and (b) modified hop-distance network topology.**

It is important to note that the hop-distance network would always result in a conservative and robust evacuation route plan, which is desirable in evacuation planning. This network would overestimate the network distance by at most 1 hop and underestimates the vehicle range at most 1 hop due to the approximation. For the case of 100 mi evacuation spacing, roads with a distance of 201 mi would be approximated as 3 hops, and vehicles with a 299 mi range would be approximated as 2 hops. Emergency planners could mitigate conservative evacuation route design by choosing a shorter controlled evacuation distance.

By applying the hop-distance network in our model, we eliminate the distance link parameter and simplify the distance measurement based on the number of hops embedded directly in the network. This hop assumption would be beneficial for emergency planners to guarantee that refueling access is satisfied within the controlled evacuation distance. In a real-world application of our model, the hop assumption is desirable and informative for emergency planners, who can make educated evacuation decisions (e.g., signage placement, routing, and fuel logistics) within a shelter's-controlled radius (FHWA, 2013). This assumption is also practical to implement in zone-based evacuation planning, where zone clusters are determined considering a controlled safe distance (Hasan and Van Hentenryck, 2020, 2021; Hsu and Peeta, 2014).

**Vehicle Driving Range Implications in Evacuation Route Planning**
The alternative fuel vehicle driving range is one of the major hindrances of adopting alternative fuel vehicles. The numerical experiments deliver two crucial findings: i) as the driving range improves, the evacuation performance improves; ii) for larger driving ranges, vehicles are less dependent on refueling infrastructure availability when evacuating. We discuss the implications of our findings under four points of view: feasibility, performance, trends, and applicability.

*Evacuation Feasibility.* If the driving range of alternative fuel vehicles is small, the evacuation plan could be infeasible. **Figure 9** shows that conventional evacuation planning becomes infeasible for alternative fuel vehicles with a low driving range.

*Evacuation Performance.* We examine the evacuation performance improvement as the driving ranges increase. First, we consider three types of alternative fuel vehicles, i.e., short-range, medium-range, and long-range vehicles. Then, we evaluate the evacuation performance improvement as we relax the driving range by 1 hop, 2 hops, and 3 hops. **Table 7** shows the improvement of the evacuation performance due to larger driving ranges for the three types of vehicles. There is a 62.71% improvement in evacuation performance as we relax the driving range of short-range vehicles by 1 hop; the same improvement drops to 19.13% for mid-range vehicles, while no improvement is observed for long-range vehicles. This indicates that driving range improvements would be significantly beneficial for the evacuation of short-range vehicles. Providing investment on improving the driving range of short-range vehicles as well as incentive programs for mid-range vehicle adoption would be beneficial in improving the region's evacuation performance.

**Table 7 The Improvement of Evacuation Performance due to Vehicle Driving Range Improvement**

| Vehicle Driving Range Bound Improvement | Evacuation Performance Improvement | | |
|---|---|---|---|
| | Short-Range Vehicles (3-5 hops) | Mid-Range Vehicles (5-8 hops) | Long-Range Vehicles (> 8 hops) |
| 1 hop | 62.71% | 19.13% | 0% |
| 2 hops | 64.03% | 32.96% | 0% |
| 3 hops | 71.58% | 41.85% | 0% |

*Range Improvement Trend.* Since the introduction of alternative fuel vehicles, their ranges have been increasing rapidly. For instance, 2021 electric vehicle (EV) models could offer a maximum range of 405 miles and a median range of 234 miles. Empirical data shows that most 2021 EVs would have a median range of approximately 60% of gasoline vehicles (US Department of Energy, 2022a). There is also a significant improvement of hydrogen fuel cell vehicles driving range. The 2022 Toyota Mirai could offer a median range of 337 miles, 20% less compared to gasoline vehicles (US Department of Energy, 2022b). This range-increasing trend is beneficial to evacuation planning and its performance. However, the range improvement feasibility often depends on automanufacturers and technological innovations, which are dynamic and uncertain.

*Long Range Applicability.* We argue that increased driving ranges of alternative fuel vehicles are not always a viable option. Large range alternative fuel vehicles are popular with early wealthy adopters but often expensive and would not be a vehicle ownership alternative for low- and middle-income communities (Lee et al., 2019). Driving shorter-range vehicles is riskier during evacuations; however, it is sufficient and socially optimal for the daily mobility needs of their drivers (Kontou et al., 2017, 2015). Shorter-range and mid-range vehicles are often promoted with governmental rebate programs (Illinois General Assembly, 2022; Texas Commission on Environmental Quality, 2022), making these vehicles more accessible for mass market adoption. Therefore, emergency planners need to understand the impact of driving range on evacuation planning and leverage our critical findings to strategically adjust the evacuation route design to address current vehicle technology barriers. Emergency planners must design unique evacuation routing plans for a diverse set of vehicle technologies that accommodates the needs of a population of diverse socio-economic characteristics.

**Refueling Infrastructure Implications in Evacuation Route Planning**

The refueling infrastructure of alternative fuel vehicles is significantly sparser compared to gasoline stations in the US (US Department of Energy, 2021a). Refueling station accessibility issues is a mobility barrier for alternative fuel vehicles during both their daily activities and emergency conditions. The numerical experiments highlight that a denser alternative refueling infrastructure improves the average evacuation performance, but the siting of stations matters. We comment on the implications of our findings under four points of view: feasibility, performance, applicability, and stakeholders.

*Evacuation Feasibility.* In **Figure 11**, we generate 32 refueling deployment cases to examine the impact of refueling infrastructure siting on evacuation planning. Our findings show that 17 of the 32 generated deployment cases with a number of 3-5 stations could result in a feasible evacuation route for alternative fuel vehicles with a minimum range of 3 hops (short-range vehicle). The remaining deployment cases, referring to 1-2 stations scenarios, could provide a feasible evacuation routing plan when the minimum range of the vehicle is 5 hops (mid-range vehicle). Denser refueling networks would significantly benefit the short-range vehicle drivers. Strategically mapping refueling stations placement would be a crucial component of supporting the evacuation process.

*Evacuation Performance.* We examine the effect of refueling station density on the evacuation performance by improving its numbers by 1, 2, 3, 4, and 5 stations. **Table 8** shows the summary of the refueling infrastructure density effect on the evacuation performance. On average, the evacuation performance improves 13.53% by adding 1 station, 36.87% by adding 3 stations, and 52.36% by adding 5 stations. Note that the installation location matters; thus, we provide the minimum and maximum range of improvement in our network experiments. Following the Electric Vehicle Charging Action Plan (2021) and California Clean Transportation Program for hydrogen-fuel stations (2022), there is an unprecedented opportunity to achieve denser refueling networks for alternative fuel vehicles in the near future, which is beneficial for emergency planning.

**Table 8 The Improvement of Evacuation Performance due to Refueling Infrastructure Investment**

| Refueling Infrastructure Investment | Evacuation Performance Improvement | | |
|---|---|---|---|
| | Average Improvement | Minimum Improvement | Maximum Improvement |
| 1 station | 13.53% | 4.89% | 20.68% |
| 2 stations | 25.11% | 10.93% | 33.86% |
| 3 stations | 36.78% | 27.96% | 46.51% |
| 4 stations | 44.53% | 39.93% | 49.13% |
| 5 stations | 52.36% | 52.36% | 52.36% |

*Infrastructure Applicability.* Our experiments point out that deploying refueling stations closer to the shelter would improve the evacuation performance. Such a placement would allow evacuees to fulfill their required refueling needs using the shorter route without a major detour. For instance, one of the numerical cases shows that adding one station closer to the shelter might improve evacuation time by 37.25%, while adding a station closer to the network center only improves it by 3.12%. The siting of refueling stations closer to the shelter would be beneficial for the evacuation operation. However, emergency planners do not often make these decisions since refueling infrastructure networks are installed and operated by industry.

*Critical Stakeholder.* Studies in the literature have contributed to planning refueling stations deployment in cities (Ghamami et al., 2015; He et al., 2015) and corridors (Erdoğan et al., 2022; Ghamami et al., 2020; Kim and Kuby, 2012; Nie and Ghamami, 2013). In practice, there are ambitious support packages from federal and state governments, investing in alternative refueling stations. For instance, DOE and DOT are investing $5 billion in nationwide EV charging infrastructure networks (US Department of Energy, 2022c). In addition, the California Energy Commission commits to allocating up to $20 million per year to fund the hydrogen fueling stations in their state (California Air Resources Board, 2021). There is a plethora of innovations on EV charging infrastructure, including wireless on-road charging (Liu and Wang, 2017; Sun et al., 2020), mobile charging (Huang et al., 2015), and peer-to-peer charging (Chakraborty et al., 2020); which could assist with emergencies and evacuations. The refueling station planning decisions do not currently consider the evacuation demand; instead, they often focus on meeting habitual refueling demand and profit maximization. While the deployment of refueling stations is often an industry-led endeavor (Gnann et al., 2018; Nicholas, 2019), public agencies often participate by subsidizing siting. Alternative refueling station investment is a multi-stakeholder decision problem, including federal, state, local governments, industry, and other entities like utilities. In 2011, regional coordination and collaboration were initiated between the Northeast and Mid-Atlantic states to develop partnerships with the private sector, utilities, and other public entities in the deployment planning of electric charging stations (Transportation and Climate Initiative, 2022). Emergency agencies must also coordinate with such entities that make decisions that impact evacuation route designs, understand how each decision affects emergency management, and denote the role they can play in refueling station investment allocation. A critical recommendation is enhanced communication and collaboration between evacuation coordination agencies, government agencies, and the alternative refueling infrastructure providers.

**Vehicle Heterogeneity Impact in Evacuation Route Planning**

Our proposed model considers the heterogeneity of vehicle fuel types. Due to the projected rise of various alternative fuel vehicle types' adoption in the future, emergency planners need to strategically manage and guarantee a safe and prompt evacuation plan for each vehicle fuel type. The vehicle heterogeneity elevates the planning complexity since emergency planners must provide and monitor different routes for each vehicle type during evacuations. Key evacuation signs placement, traffic control, evacuation scheduling, and emergency logistics management for each vehicle fuel type are crucial in supporting such evacuations.

Furthermore, the vehicle heterogeneity could be extended into the varieties of vehicle ranges, such as short-range, mid-range, and long-range vehicles for the same vehicle fuel type. Our proposed model can easily handle this extension by defining the $|K|$-evacuation tree routes problem for $|K|$ predetermined driving range bounds. Our proposed approach allows evacuation planners and managers to determine the diverse vehicle fuel types they want to incorporate in their plans.

Consequently, the heterogeneity of alternative fuel vehicles would lead to further research questions, such as what the important classifications are, how to determine these classes, who determines the classifications, and how many classes need to be considered in evacuation planning. Our proposed evacuation model contributes to initiating discussion on how vehicle heterogeneity impacts preemptive evacuation planning. The proposed model can also serve as a platform of enumeration for deciding the proper heterogeneous vehicle classes for evacuation planning.

## CONCLUSIONS

In summary, the alternative fuel vehicles adoption has grown rapidly in the US (**US Department of Energy, 2021e**) and was followed by increased governmental support (**State of Illinois, 2021; The White House, 2021**) that aims to further incentivize alternative fuel vehicles ownership. Their advent raises concerns on how such emerging vehicle technologies will impact evacuation route planning (**Adderly et al., 2018; Feng et al., 2020**). To date, there are very few studies that focus on coupling evacuation route planning and alternative fuel vehicles, but our work is the first aiming to design preemptive evacuation routes for a diverse set of vehicle fuel types.

This paper has presented a novel evacuation route planning problem for alternative fuel vehicles, by developing a $|K|$-evacuation tree routes model (rooted at the safety node) for a set of different vehicle fuel types $K$, as well as a solution algorithm based on matheuristic branch-and-price. The model incorporates driving range constraints and refueling needs to enable alternative fuel vehicles to reach safety in a reliable manner. The mathematical formulation minimizes the total evacuation time, which is the summation of the total travel time and total refueling time. In the proposed formulation, we represent the vehicles' driving range as the allowable distance (in number of hops) from each origin node to the safety directly without refueling. We impose hop constraints to capture the refueling needs of each alternative fuel vehicle type. We present numerical examples where we apply the proposed evacuation planning problem to the Sioux Falls and South Florida transportation networks with existing alternative refueling infrastructure deployment.

The experiments show that the evacuation routing plans for each vehicle fuel type naturally change as we consider more heterogeneous driving ranges and different refueling needs on the path to reach safety. The evacuation route for each vehicle fuel type is different and could prove infeasible to other vehicle fuel types due to different availability and placement of their refueling infrastructure. Besides, the experiments show that the driving range constraint plays a vital role in planning preemptive evacuation routes. The total evacuation time of the vehicles in the network decreases as the driving range increases and the availability and placement of refueling stations of alternative vehicle fuel types are important in designing faster evacuations. Our experiments indicate that denser deployment of charging infrastructure would support faster and more reliable evacuation planning, when deployed at appropriate network nodes, but the siting of refueling station matters.

The advantage of the proposed evacuation planning problem and its solutions approach is that it considers driving range constraints and refueling needs of alternative fuel vehicles while providing a seamless, reliable, and feasible evacuation plan for multiple alternative vehicle fuel types. Limitations that need to be addressed in future studies include accounting for the heterogeneity of the remaining fuel/battery range of each vehicle type during the evacuation. In this study, we assume that evacuees originating from every node that is located further from the shelter than the allowable driving range need to recharge before reaching safety. However, this assumption may not necessarily hold given that vehicles may have sufficient remaining driving range to reach safety regardless of their current position. In this study, we assume that the charging time is proportional to a constant rate and the distance from each origin node to the safety node, without any queueing time and delays considerations. It is implied that every refueling station has an unlimited capacity to serve incoming vehicles. However, alternative fuel infrastructure might have a limited number of refueling hoses or charging ports to accommodate evacuees' vehicles, which could result in queueing at a refueling infrastructure node. Future modeling extensions include dynamic models that capture such delays. Additionally, we realize the need to test and implement the proposed formulation to real-world large-scale transportation networks and realistic evacuation scenario. This is a recommended future direction for expanding contributions in this field.

## ACKNOWLEDGMENTS
This research did not receive funding from agencies in the public, commercial, or not-for-profit sectors. The authors thank Dr. Xiaotong Sun and Dr. Zhibin Chen for providing the data on the South Florida transportation network. The authors confirm full personal access to all aspects of the writing and research process and take full responsibility for the paper.

## REFERENCES

Achrekar, O., Vogiatzis, C., 2018. Evacuation trees with contraflow and divergence considerations, in: Dynamics of Disasters. Springer International Publishing, pp. 1–46. https://doi.org/10.1007/978-3-319-97442-2_1

Adderly, S.A., Manukian, D., Sullivan, T.D., Son, M., 2018. Electric vehicles and natural disaster policy implications. Energy Policy 112, 437–448. https://doi.org/10.1016/j.enpol.2017.09.030

Agrawal, S., Zheng, H., Peeta, S., Kumar, A., 2016. Routing aspects of electric vehicle drivers and their effects on network performance. Transportation Research Part D: Transport and Environment 46, 246–266. https://doi.org/10.1016/j.trd.2016.04.002

Andreas, A.K., Smith, J.C., 2009. Decomposition Algorithms for the Design of a Nonsimultaneous Capacitated Evacuation Tree Network. Networks 53, 91–103. https://doi.org/10.1002/net.20278

Apte, A., 2009. Humanitarian logistics: A new field of research and action. Foundations and Trends in Technology, Information and Operations Management 3, 1–100. https://doi.org/10.1561/0200000014

Archetti, C., Bianchessi, N., Speranza, M.G., 2011. A column generation approach for the split delivery vehicle routing problem. Networks 58, 241–254. https://doi.org/10.1002/NET.20467

Archetti, C., Speranza, M.G., 2014. A survey on matheuristics for routing problems. EURO Journal on Computational Optimization 2, 223–246. https://doi.org/10.1007/S13675-014-0030-7/FIGURES/1

Barnhart, C., Johnson, E.L., Nemhauser, G.L., Savelsbergh, M.W.P., Vance, P.H., 1998. Branch-and-price: Column generation for solving huge integer programs. Operations Research 46, 316–329. https://doi.org/10.1287/OPRE.46.3.316

Bayram, V., 2016. Optimization models for large scale network evacuation planning and management: A literature review. Surveys in Operations Research and Management Science 21, 63–84. https://doi.org/10.1016/j.sorms.2016.11.001

Bayram, V., Tansel, B.T., Yaman, H., 2015. Compromising system and user interests in shelter location and evacuation planning. Transportation Research Part B: Methodological 72, 146–163. https://doi.org/10.1016/j.trb.2014.11.010

Bender, M.A., Knutson, T.R., Tuleya, R.E., Sirutis, J.J., Vecchi, G.A., Garner, S.T., Held, I.M., 2010. Modeled impact of anthropogenic warming on the frequency of intense Atlantic hurricanes. Science (1979) 327, 454–458. https://doi.org/10.1126/science.1180568

California Air Resources Board, 2021. CEC and CARB Joint Agency Report on Assembly Bill 8.

California Energy Commission, 2021. Zero Emission Vehicle and Infrastructure Statistics [WWW Document]. URL https://www.energy.ca.gov/data-reports/energy-insights/zero-emission-vehicle-and-charger-statistics (accessed 7.25.21).

California FPD, 2022. Fire Evacuation Plans | San Ramon Valley Fire Protection District [WWW Document]. URL https://www.firedepartment.org/our-district/fire-and-life-safety/fire-evacuation-plans (accessed 2.26.22).

California Legislature, 2022. Bill Text - AB-2562 Clean Transportation Program: hydrogen-fueling stations. [WWW Document]. URL https://leginfo.legislature.ca.gov/faces/billTextClient.xhtml?bill_id=202120220AB2562 (accessed 2.23.22).

Campos, V., Bandeira, R., Bandeira, A., 2012. A Method for Evacuation Route Planning in Disaster Situations. Procedia - Social and Behavioral Sciences 54, 503–512. https://doi.org/10.1016/J.SBSPRO.2012.09.768

Chakraborty, P., Parker, R., Hoque, T., Cruz, J., Bhunia, S., 2020. P2C2: Peer-to-Peer Car Charging. IEEE Vehicular Technology Conference 2020-May. https://doi.org/10.1109/VTC2020-SPRING48590.2020.9128955

Chen, X., Kwan, M.P., Li, Q., Chen, J., 2012. A model for evacuation risk assessment with consideration of pre- and post-disaster factors. Computers, Environment and Urban Systems 36, 207–217. https://doi.org/10.1016/j.compenvurbsys.2011.11.002

Chen, X., Meaker, J.W., Zhan, F.B., 2006. Agent-based modeling and analysis of hurricane evacuation procedures for the Florida Keys. Natural Hazards 38, 321–338. https://doi.org/10.1007/s11069-005-0263-0

Chen, X., Zhan, F.B., 2014. Agent-based modeling and simulation of urban evacuation: relative effectiveness of simultaneous and staged evacuation strategies, in: Agent-Based Modeling and Simulation. Palgrave Macmillan UK, pp. 78–96. https://doi.org/10.1057/9781137453648_6

Cova, T.J., Johnson, J.P., 2003. A network flow model for lane-based evacuation routing. Transportation Research Part A: Policy and Practice 37, 579–604. https://doi.org/10.1016/S0965-8564(03)00007-7

Davis, S.C., Boundy, R.G., 2021. Transportation Energy Data Book: Edition 39.

Desaulniers, G., Desrosiers, J., Solomon, M.M., 2005. Column generation, Column Generation. Springer US. https://doi.org/10.1007/B135457

Desrosiers, J., Lübbecke, M.E., 2005. A primer in column generation. Column Generation 1–32. https://doi.org/10.1007/0-387-25486-2_1

Ebihara, M., Ohtsuki, A., Iwaki, H., 1992. A Model for Simulating Human Behavior During Emergency Evacuation Based on Classificatory Reasoning and Certainty Value Handling. Computer‐Aided Civil and Infrastructure Engineering 7, 63–71. https://doi.org/10.1111/j.1467-8667.1992.tb00417.x

Erdoğan, S., Çapar, İsmail, Çapar, İbrahim, Nejad, M.M., 2022. Establishing a statewide electric vehicle charging station network in Maryland: A corridor-based station location problem. Socio-Economic Planning Sciences 79. https://doi.org/10.1016/j.seps.2021.101127

Erdoğan, S., Miller-Hooks, E., 2012. A Green Vehicle Routing Problem. Transportation Research Part E: Logistics and Transportation Review 48, 100–114. https://doi.org/10.1016/j.tre.2011.08.001

Eshghi, K., Larson, R.C., 2008. Disasters: Lessons from the past 105 years. Disaster Prevention and Management: An International Journal 17, 62–82. https://doi.org/10.1108/09653560810855883

FDH, 2022. Emergency Preparedness and Response | Florida Department of Health in St Johns [WWW Document]. URL https://stjohns.floridahealth.gov/programs-and-services/emergency-preparedness-and-response/ (accessed 2.26.22).

FDOT, 2021. Emergency Shoulder Use (ESU) [WWW Document]. URL https://www.fdot.gov/emergencymanagement/esu/default.shtm (accessed 2.26.22).

Federal Highway Administration, 2020. Good Practices in Transportation Evacuation Preparedness and Response - Phase 1 - Preparation and Activation - FHWA Emergency Transportation Operations [WWW Document]. Emergency Transportation Operations. URL https://ops.fhwa.dot.gov/publications/fhwahop09040/phase1.htm (accessed 2.28.22).

Federal Highway Administration, 2017. Evacuation Stakeholders' Roles and Responsibilities [WWW Document]. outes to Effective Evacuation Planning Primer Series: Using Highways During Evacuation Operations for Events with Advance Notice. URL https://ops.fhwa.dot.gov/publications/evac_primer/05_evacuation.htm (accessed 2.28.22).

Federal Highway Administration, 2006. Assessment of the State of the Practice and State of the Art in Evacuation Transportation Management.

FEMA, 2022. Declared Disasters | FEMA.gov [WWW Document]. URL https://www.fema.gov/disaster/declarations (accessed 2.24.22).

FEMA, 2019. Planning Considerations: Evacuation and Shelter-In-Place. Federal Emergency Management Agency.

Feng, K., Lin, N., Xian, S., Chester, M. v., 2020. Can we evacuate from hurricanes with electric vehicles? Transportation Research Part D: Transport and Environment 86, 102458. https://doi.org/10.1016/j.trd.2020.102458

FHWA, 2013. Chapter 2I. Emergency Managament Signing [WWW Document]. Manual on Uniform Traffic Control Devices. URL https://mutcd.fhwa.dot.gov/htm/2003r1/part2/part2i.htm (accessed 2.21.22).

Fischetti, M., Lancia, G., Serafini, P., 2002. Exact Algorithms for Minimum Routing Cost Trees. Networks 39, 161–173. https://doi.org/10.1002/net.10022

Franke, T., Krems, J.F., 2013. Understanding charging behaviour of electric vehicle users. Transportation Research Part F: Traffic Psychology and Behaviour 21, 75–89. https://doi.org/10.1016/j.trf.2013.09.002

Gao, Y., Chiu, Y.C., Wang, S., Küçükyavuz, S., 2010. Optimal refueling station location and supply planning for hurricane evacuation. Transportation Research Record 56–64. https://doi.org/10.3141/2196-06

Ghamami, M., Kavianipour, M., Zockaie, A., Hohnstadt, L.R., Ouyang, Y., 2020. Refueling infrastructure planning in intercity networks considering route choice and travel time delay for mixed fleet of electric and conventional vehicles. Transportation Research Part C: Emerging Technologies 120, 102802. https://doi.org/10.1016/j.trc.2020.102802

Ghamami, M., Nie, Y., Zockaie, A., 2015. Planing Charging Infrastructure for Plug-in Electric Vehicles in City Centers. International Journal of Sustainable Transportation.

https://doi.org/http://www.tandfonline.com/doi/abs/10.1080/15568318.2014.937840#.VjJIkRFdEuU

Gnann, T., Funke, S., Jakobsson, N., Plötz, P., Sprei, F., Bennehag, A., 2018. Fast charging infrastructure for electric vehicles: Today's situation and future needs. Transportation Research Part D: Transport and Environment 62, 314–329. https://doi.org/10.1016/J.TRD.2018.03.004

Gouveia, L., Paias, A., Sharma, D., 2008. Modeling and solving the rooted distance-constrained minimum spanning tree problem. Computers and Operations Research 35, 600–613. https://doi.org/10.1016/j.cor.2006.03.022

Gouveia, L., Simonetti, L., Uchoa, E., 2011. Modeling hop-constrained and diameter-constrained minimum spanning tree problems as Steiner tree problems over layered graphs. Mathematical Programming 128, 123–148. https://doi.org/10.1007/s10107-009-0297-2

Gurobi Optimization, LLC, 2021. Gurobi Optimizer Reference Manual.

Hagberg, A.A., Schult, D.A., Swart, P.J., 2008. Exploring Network Structure, Dynamics, and Function using NetworkX, in: Varoquaux, G., Vaught, T., Millman, J. (Eds.), Proceedings of the 7th Python in Science Conference. Pasadena, CA USA, pp. 11–15.

Hamacher, H.W., Tjandra, S.A., 2002. Mathematical modelling of evacuation problems: a state of the art. Pedestrian and Evacuation Dynamics 24, 227–266. https://doi.org/citeulike-article-id:6650160

Hamacher, H.W., Tjandra, S.A., 2001. Mathematical Modelling of Evacuation Problems: A State of Art, Berichte des Fraunhofer ITWM, Nr.

Hasan, M.H., Van Hentenryck, P., 2020. Large-scale zone-based evacuation planning-Part I: Models and algorithms. https://doi.org/10.1002/net.21981

Hasan, Mohd.H., Van Hentenryck, P., 2021. Large-scale zone-based evacuation planning, Part II: Macroscopic and microscopic evaluations. Networks 77, 341–358. https://doi.org/10.1002/NET.21980

He, F., Wu, D., Yin, Y., Guan, Y., 2013. Optimal deployment of public charging stations for plug-in hybrid electric vehicles. Transportation Research Part B: Methodological 47, 87–101. https://doi.org/10.1016/j.trb.2012.09.007

He, F., Yin, Y., Lawphongpanich, S., 2014. Network equilibrium models with battery electric vehicles. Transportation Research Part B: Methodological 67, 306–319. https://doi.org/10.1016/j.trb.2014.05.010

He, F., Yin, Y., Zhou, J., 2015. Deploying public charging stations for electric vehicles on urban road networks. Transportation Research Part C: Emerging Technologies 60, 227–240. https://doi.org/10.1016/J.TRC.2015.08.018

He, X., Zhang, S., Wu, Y., Wallington, T.J., Lu, X., Tamor, M.A., McElroy, M.B., Zhang, K.M., Nielsen, C.P., Hao, J., 2019. Economic and Climate Benefits of Electric Vehicles in China, the United States, and Germany. Environmental Science and Technology 53, 11013–11022. https://doi.org/10.1021/acs.est.9b00531

Hori, M., Schafer, M.J., Bowman, D.J., 2009. Displacement dynamics in Southern Louisiana after Hurricanes Katrina and Rita. Population Research and Policy Review 28, 45–65. https://doi.org/10.1007/s11113-008-9118-1

Hsu, Y.T., Peeta, S., 2014. Risk-based spatial zone determination problem for stage-based evacuation operations. Transportation Research Part C: Emerging Technologies 41, 73–89. https://doi.org/10.1016/J.TRC.2014.01.013

Huang, S., He, L., Gu, Y., Wood, K., Benjaafar, S., 2015. Design of a Mobile Charging Service for Electric Vehicles in an Urban Environment. IEEE Transactions on Intelligent Transportation Systems 16, 787–798. https://doi.org/10.1109/TITS.2014.2341695

Illinois General Assembly, 2022. 415 ILCS 120 - Electric Vehicle Rebate Act. [WWW Document]. URL https://www.ilga.gov/legislation/ilcs/ilcs3.asp?ActID=1608&ChapterID=36 (accessed 2.21.22).

Jacobson, M.Z., 2009. Review of solutions to global warming, air pollution, and energy security. Energy and Environmental Science 2, 148–173. https://doi.org/10.1039/b809990c

Jiang, N., Xie, C., Waller, S., 2012. Path-constrained traffic assignment. Transportation Research Record 25–33. https://doi.org/10.3141/2283-03

Kim, J.G., Kuby, M., 2012. The deviation-flow refueling location model for optimizing a network of refueling stations. International Journal of Hydrogen Energy 37, 5406–5420. https://doi.org/10.1016/j.ijhydene.2011.08.108

Kim, S., Shekhar, S., 2005. Contraflow network reconfiguration for evacuation planning: A summary of results. GIS: Proceedings of the ACM International Symposium on Advances in Geographic Information Systems 250–259.

Kim, S., Shekhar, S., Min, M., 2008. Contraflow transportation network reconfiguration for evacuation route planning. IEEE Transactions on Knowledge and Data Engineering 20, 1115–1129. https://doi.org/10.1109/TKDE.2007.190722

Kontou, E., Yin, Y., Lin, Z., 2015. Socially optimal electric driving range of plug-in hybrid electric vehicles. Transportation Research Part D: Transport and Environment 39, 114–125. https://doi.org/10.1016/j.trd.2015.07.002

Kontou, E., Yin, Y., Lin, Z., He, F., 2017. Socially Optimal Replacement of Conventional with Electric Vehicles for the U.S. Household Fleet. International Journal of Sustainable Transportation 11, 749–763. https://doi.org/10.1080/15568318.2017.1313341

Kuby, M., Lim, S., 2005. The flow-refueling location problem for alternative-fuel vehicles. Socio-Economic Planning Sciences 39, 125–145. https://doi.org/10.1016/j.seps.2004.03.001

LCG Traffic and Transportation Department, 2022. Lafayette Hurricane Evacuation Routes [WWW Document]. URL http://lafayetteohsep.org/Pages/Lafayette-Hurricane-Evacuation-Routes (accessed 2.26.22).

Lee, J.H., Hardman, S.J., Tal, G., 2019. Who is buying electric vehicles in California? Characterising early adopter heterogeneity and forecasting market diffusion. Energy Research & Social Science 55, 218–226. https://doi.org/10.1016/J.ERSS.2019.05.011

Li, C., Cao, Y., Zhang, M., Wang, J., Liu, J., Shi, H., Geng, Y., 2015. Hidden benefits of electric vehicles for addressing climate change. Scientific Reports 5, 8–11. https://doi.org/10.1038/srep09213

Lindell, M.K., Murray-tuite, P., Wolshon, B., Baker, E.J., 2019. Large-Scale Evacuation : The Analysis, Modeling, and Management of Emergency Relocation from Hazardous Areas Michael. Routledge, New York.

Lindell, M.K., Prater, C.S., 2007. Critical behavioral assumptions in evacuation time estimate analysis for private vehicles: Examples from hurricane research and planning. Journal of Urban Planning and Development 133, 18–29. https://doi.org/10.1061/(ASCE)0733-9488(2007)133:1(18)

Liu, H., Wang, D.Z.W., 2017. Locating multiple types of charging facilities for battery electric vehicles. Transportation Research Part B: Methodological 103, 30–55. https://doi.org/10.1016/J.TRB.2017.01.005

Lutsey, N., 2018. California's continued electric vehicle market development [WWW Document]. ICCT briefing. URL https://theicct.org/publication/californias-continued-electric-vehicle-market-development/

Maas, P., Iyer, S., Gros, A., Park, W., Mcgorman, L., Nayak, C., Dow, A.P., 2019. Facebook Disaster Maps: Aggregate Insights for Crisis Response & Recovery, in: Proceedings of the 16th ISCRAM Conference.

Metaxa-Kakavouli, D., Maas, P., Aldrich, D.P., 2018. How social ties influence hurricane evacuation behavior. Proceedings of the ACM on Human-Computer Interaction 2. https://doi.org/10.1145/3274391

Miotti, M., Supran, G.J., Kim, E.J., Trancik, J.E., 2016. Personal Vehicles Evaluated against Climate Change Mitigation Targets. Environmental Science and Technology 50, 10795–10804. https://doi.org/10.1021/acs.est.6b00177

Murray-Tuite, P., Wolshon, B., 2013. Evacuation transportation modeling: An overview of research, development, and practice. Transportation Research Part C: Emerging Technologies 27, 25–45. https://doi.org/10.1016/j.trc.2012.11.005

National Oceanic and Atmospheric Administration, 2020. Tropical Cyclone Climatology [WWW Document]. URL https://www.nhc.noaa.gov/climo/ (accessed 2.24.22).

Ng, M.W., Park, J., Waller, S.T., 2010. A Hybrid Bilevel Model for the Optimal Shelter Assignment in Emergency Evacuations. Computer-Aided Civil and Infrastructure Engineering 25, 547–556. https://doi.org/10.1111/j.1467-8667.2010.00669.x

Ng, M.W., Waller, S.T., 2009. The evacuation optimal network design problem: Model formulation and comparisons. Transportation Letters 1, 111–119. https://doi.org/10.3328/TL.2009.01.02.111-119

Nicholas, M., 2019. Estimating electric vehicle charging infrastructure costs across major U.S. metropolitan areas.

Nie, Y., Ghamami, M., 2013. A corridor-centric approach to planning electric vehicle charging infrastructure. Transportation Research Part B: Methodological 57, 172–190. https://doi.org/10.1016/j.trb.2013.08.010

NYC Office of Climate and Sustainability, 2021. New York City's Net-Zero Carbon Target for 2050 Is Achievable, Study Finds - Sustainability [WWW Document]. URL https://www1.nyc.gov/site/sustainability/our-programs/carbon-neutral-nyc-pr-04-15-2021.page (accessed 6.25.21).

Office of Energy Efficiency & Renewable Energy, 2018. FOTW #1010: All-Electric Light Vehicle Ranges Can Exceed Those of Some Gasoline Light Vehicles [WWW Document]. Vehicle Technologies Office. URL https://www.energy.gov/eere/vehicles/articles/fotw-1010-january-1-2018-all-electric-light-vehicle-ranges-can-exceed-those (accessed 6.25.21).

Ogier, R.G., 1988. Minimum-delay routing in continuous-time dynamic networks with Piecewise-constant capacities. Networks 18, 303–318. https://doi.org/10.1002/net.3230180405

Opasanon, S., Miller-Hooks, E., 2010. Noisy genetic algorithm for stochastic, time-varying minimum time network flow problem. Transportation Research Record 75–82. https://doi.org/10.3141/2196-08

Pillac, V., Van Hentenryck, P., Even, C., 2016. A conflict-based path-generation heuristic for evacuation planning. Transportation Research Part B: Methodological 83, 136–150. https://doi.org/10.1016/j.trb.2015.09.008

Rezvani, Z., Jansson, J., Bodin, J., 2015. Advances in consumer electric vehicle adoption research: A review and research agenda. Transportation Research Part D: Transport and Environment 34, 122–136. https://doi.org/10.1016/j.trd.2014.10.010

Roberts, D., 2019. New York passes the country's most ambitious climate target [WWW Document]. Vox. URL https://www.vox.com/energy-and-environment/2019/6/20/18691058/new-york-green-new-deal-climate-change-cuomo (accessed 6.25.21).

Sabbaghtorkan, M., Batta, R., He, Q., 2022. On the analysis of an idealized model to manage gasoline supplies in a short-notice hurricane evacuation. OR Spectrum. https://doi.org/10.1007/s00291-022-00665-0

Sarma, D., Das, A., Bera, U.K., 2020. Uncertain demand estimation with optimization of time and cost using Facebook disaster map in emergency relief operation. Applied Soft Computing 87, 105992. https://doi.org/10.1016/J.ASOC.2019.105992

Sbayti, H., Mahmassani, H.S., 2006. Optimal Scheduling of Evacuation Operations. Transportation Research Record: Journal of the Transportation Research Board 1964, 238–246. https://doi.org/10.1177/0361198106196400126

Sheffi, Y., 1985. Urban Transportation Networks. Prentice-Hall, Englewood Cliffs, New Jersey, USA.

State of Illinois, 2021. Illinois Press Release: Gov. Pritzker Signs Transformative Legislation Establishing Illinois as a National Leader on Climate Action [WWW Document]. Illinois Press Release. URL https://www.illinois.gov/news/press-release.23893.html (accessed 2.28.22).

Sun, X., Chen, Z., Yin, Y., 2020. Integrated planning of static and dynamic charging infrastructure for electric vehicles. Transportation Research Part D: Transport and Environment 83, 102331. https://doi.org/10.1016/j.trd.2020.102331

Tessum, C.W., Hill, J.D., Marshall, J.D., 2014. Life cycle air quality impacts of conventional and alternative light-duty transportation in the United States. Proc Natl Acad Sci U S A 111, 18490–18495. https://doi.org/10.1073/pnas.1406853111

Texas Commission on Environmental Quality, 2022. Light-Duty Motor Vehicle Purchase or Lease Incentive Program [WWW Document]. URL https://www.tceq.texas.gov/airquality/terp/ld.html (accessed 2.21.22).

The White House, 2021. FACT SHEET: President Biden Announces Steps to Drive American Leadership Forward on Clean Cars and Trucks | The White House [WWW Document]. The White House Briefing Room: Statements and Releases. URL https://www.whitehouse.gov/briefing-room/statements-releases/2021/08/05/fact-sheet-president-biden-announces-steps-to-drive-american-leadership-forward-on-clean-cars-and-trucks/ (accessed 2.28.22).

Tilk, C., Irnich, S., 2018. Combined column-and-row-generation for the optimal communication spanning tree problem. Computers and Operations Research 93, 113–122. https://doi.org/10.1016/j.cor.2018.01.003

Townsend, F.F., 2006. Hurricane Katrina: Lessons Learned - Letter from Frances Fragos Townsend [WWW Document]. URL https://georgewbush-whitehouse.archives.gov/reports/katrina-lessons-learned/letter.html (accessed 2.26.22).

Transportation and Climate Initiative, 2022. Northeast Electric Vehicle Network in Action [WWW Document]. URL https://www.transportationandclimate.org/content/northeast-electric-vehicle-network (accessed 5.19.22).

Transportation Networks for Research Core Team, 2021. Transportation Networks for Research [WWW Document]. URL https://github.com/bstabler/TransportationNetworks/tree/master/SiouxFalls (accessed 6.27.21).

TxDOT, 2021. I-10 Contraflow.

US Department of Energy, 2022a. FOTW #1221, January 17, 2022 [WWW Document]. Vehicle Technology Offices. URL https://www.energy.gov/eere/vehicles/articles/fotw-1221-january-17-2022-model-year-2021-all-electric-vehicles-had-median (accessed 2.21.22).

US Department of Energy, 2022b. Compare Fuel Cell Vehicles [WWW Document]. Fuel Economy. URL https://www.fueleconomy.gov/feg/fcv_sbs.shtml (accessed 2.21.22).

US Department of Energy, 2022c. President Biden, DOE and DOT Announce $5 Billion over Five Years for National EV Charging Network [WWW Document].

US Department of Energy, 2021a. Electric Vehicle Charging Station Locations [WWW Document]. Alternative Fuels Data Center. URL https://afdc.energy.gov/fuels/electricity_locations.html#/find/nearest?fuel=ELEC (accessed 6.25.21).

US Department of Energy, 2021b. Fuel Cell Electric Vehicles [WWW Document]. Alternative Fuels Data Center. URL https://afdc.energy.gov/vehicles/fuel_cell.html (accessed 7.13.21).

US Department of Energy, 2021c. Compressed Natural Gas Fueling Stations [WWW Document]. Alternative Fuels Data Center. URL https://afdc.energy.gov/fuels/natural_gas_cng_stations.html (accessed 7.13.21).

US Department of Energy, 2021d. Developing Infrastructure to Charge Plug-In Electric Vehicles [WWW Document]. Alternative Fuels Data Center. URL https://afdc.energy.gov/fuels/electricity_infrastructure.html (accessed 7.13.21).

US Department of Energy, 2021e. New Plug-in Electric Vehicle Sales in the United States Nearly Doubled from 2020 to 2021 | Department of Energy [WWW Document]. URL https://www.energy.gov/energysaver/articles/new-plug-electric-vehicle-sales-united-states-nearly-doubled-2020-2021 (accessed 3.2.22).

US Energy Information Administration, 2020. Annual Energy Outlook 2020 with projections to 2050 [WWW Document].

Vanderbeck, F., 2005. Implementing mixed integer column generation. Column Generation 331–358. https://doi.org/10.1007/0-387-25486-2_12

Vanderbeck, F., 2000. On Dantzig-Wolfe decomposition in integer programming and ways to perform branching in a branch-and-price algorithm. Operations Research 48, 111–128. https://doi.org/10.1287/OPRE.48.1.111.12453

Vogiatzis, C., Walteros, J.L., Pardalos, P.M., 2013. Evacuation Through Clustering Techniques, in: Goldengorin, B., Kalyagin, V.A., Pardalos, P.M. (Eds.), Models, Algorithms, and Technologies for Network Analysis: Proceedings of the First International Conference on Network Analysis, Springer Proceedings in Mathematics & Statistics. Springer New York, New York, NY, pp. 185–198. https://doi.org/10.1007/978-1-4614-5574-5

Wei, W., Wu, L., Wang, J., Mei, S., 2018. Network equilibrium of coupled transportation and power distribution systems. IEEE Transactions on Smart Grid 9, 6764–6779. https://doi.org/10.1109/TSG.2017.2723016

White House, 2021. FACT SHEET: The Biden-Harris Electric Vehicle Charging Action Plan | The White House [WWW Document]. URL https://www.whitehouse.gov/briefing-room/statements-releases/2021/12/13/fact-sheet-the-biden-harris-electric-vehicle-charging-action-plan/ (accessed 2.23.22).

Wolshon, B.B., 2001. "'ONE-WAY-OUT'": CONTRAFLOW FREEWAY OPERATION FOR HURRICANE EVACUATION. Natural Hazards Review 2, 105–112.

Xie, C., Lin, D.Y., Travis Waller, S., 2010. A dynamic evacuation network optimization problem with lane reversal and crossing elimination strategies. Transportation Research Part E: Logistics and Transportation Review 46, 295–316. https://doi.org/10.1016/j.tre.2009.11.004

Xie, C., Turnquist, M.A., 2011. Lane-based evacuation network optimization: An integrated Lagrangian relaxation and tabu search approach. Transportation Research Part C: Emerging Technologies 19, 40–63. https://doi.org/10.1016/j.trc.2010.03.007

Yamada, T., 1996. A network flow approach to a city emergency evacuation planning. International Journal of Systems Science 27, 931–936. https://doi.org/10.1080/00207729608929296